\documentclass[a4paper, 11 pt , english]{article}
\usepackage{ulem}
\usepackage[english]{babel}
\usepackage{amsthm}
\usepackage{pstricks,pst-node,pst-tree}
\usepackage{amssymb}
\usepackage{enumitem}
\usepackage[utf8]{inputenc}
\usepackage[T1]{fontenc} 
\usepackage{fancybox} 
\usepackage{alltt}
\usepackage{graphicx}
\usepackage{bm}
\usepackage{lmodern}           
\usepackage{babel}
\usepackage[top = 3cm, bottom = 3 cm, right = 3cm, left = 3 cm]{geometry}
\usepackage{appendix}
\usepackage{float}
\usepackage{pgf,tikz,pgfplots}
\pgfplotsset{compat=1.15}
\usepackage{mathrsfs}
\usepackage{dsfont}
\usetikzlibrary{arrows}
\usepackage{graphicx}
\usepackage{amsmath}
\usepackage{amssymb}
\usepackage{mathtools}       
\usepackage[T1]{fontenc}   
\usepackage{lmodern}     
\usepackage[utf8]{inputenc}   
\usepackage{babel}
\usepackage[linesnumbered, french]{algorithm2e}
\usepackage{physics}
\usepackage{caption}
\usepackage{subcaption} 
\usepackage{caption}
\usepackage{mathtools}
\usepackage{empheq}
\usepackage[most]{tcolorbox}
\usepackage{stmaryrd}
\usepackage{graphicx}
\usepackage{bigints}
\newtheorem{Theorem}{Theorem}[section]
\newtheorem{Definition}[Theorem]{Definition}
\newtheorem{Proposition}[Theorem]{Proposition}
\newtheorem{Assumption}[Theorem]{Assumption}
\newtheorem{Lemma}[Theorem]{Lemma}

\newtheorem{Remark}[Theorem]{Remark}

\newtheorem*{AssumptionC}{Assumption $\textbf{\rm C}_{\mathbb{W}^{\rm ad}_p}$}
\newtheorem*{AssumptionD}{Assumption $\textbf{\rm D} $}

\def \Prod{\displaystyle\prod}

\def \be{\begin{eqnarray}}
\def \ee{\end{eqnarray}}
\def \b*{\begin{align*}}
\def \e*{\end{align*}}

\def \E{\mathbb{E}}

\def \N{\mathbb{N}}
\def \P{\mathbb{P}}

\def \R{\mathbb{R}}

\def \W{\mathbb{W}}
\def \X{\mathbb{X}}

\def \[{[\,[}
\def \]{]\,]}

\def \1{{\bf 1}}

\def \ep{\hbox{}\hfill$\Box$}

\def\Bc{{\cal B}}

\def\Ec{{\cal E}}
\def\Fc{{\cal F}}
\def\Gc{{\cal G}}

\def\Lc{{\cal L}}

\def\Nc{{\cal N}}

\def\Pc{{\cal P}}

\def\Tc{{\cal T}}
\def\Uc{{\cal U}}

\def\Yc{{\cal Y}}
\def\Zc{{\cal Z}}

\def\namedlabel#1#2{\begingroup
  #2%
  \def\@currentlabel{#2}%
  \phantomsection\label{#1}\endgroup
}

\def \Lebone{{\rm Leb}_1}
\def \Lebtwo{{\rm Leb}_2}

\def \lowGmartcoupad{\underline{G}^{{\rm M}, {\rm m}}_{\W_p^{\rm ad}}}
\def \upGmartcoupad{\overline{G}^{{\rm M}, {\rm m}}_{\W_p^{\rm ad}}}

\def \lowGgenconsddist{\underline{G}^{\mathcal{C}}}
\def \upGgenconsddist{\overline{G}^{\mathcal{C}}}

\def \opedistmartcoup{\mathcal{L}_{\rm M, m}}

\def \Bad{B_{\W_p^{\rm ad}}}
\def \BadMm{B^{{\rm M}, {\rm m}}_{\W_p^{\rm ad}}}

\def \Bdisgen{B_{\rm ad}^{\varphi, \Ec}}

\def \Econsgen{\Ec}
\def \clEconsdist{\overline{\Ec}^{\rm ad}}

\def \Nad{\mathbf{N}_{\rm ad}}

\def \x{\times}

\newcommand{\vertiii}[1]{{\left\vert\kern-0.25ex\left\vert\kern-0.25ex\left\vert #1 
  \right\vert\kern-0.25ex\right\vert\kern-0.25ex\right\vert}}

\usepackage{epsfig}
\usepackage{fullpage}
\usepackage{fancyhdr}
\usepackage{moreverb}
\usepackage{xspace}
\usepackage[colorlinks,hyperindex,bookmarks,linkcolor=blue,citecolor=blue,urlcolor=blue]{hyperref}

\usepackage{wrapfig}
\usepackage{epsf}
\usepackage[maxnames = 50]{biblatex}

\title{\bf Model Risk Static-Hedging a Constrained Distributionally Robust Optimization approach}
\author{            Nathan Sauldubois\thanks{
            New York University, Tandon School of Engineering
                                ns6982@nyu.edu
                                } 
}
\date{\today}  

\DefineBibliographyStrings{english}{
  and = {\&}}

\usepackage{xcolor}
\allowdisplaybreaks
\begin{document}

\maketitle

\begin{abstract}
We investigate model risk and Distributionally Robust optimization (DRO) under marginal and martingale constraints. 
This article naturally continues the work of \citeauthor{Touzisauldubois2024ordermartingalemodelrisk} \cite{Touzisauldubois2024ordermartingalemodelrisk}, where we left open the case of static hedging with second‑period maturity vanilla options and hedging strategies involving a vanilla payoff. 
We extend the results of \citeauthor{bartl2021sensitivity} \cite{bartl2021sensitivity} and \citeauthor{bartlsensitivityadapted} \cite{bartlsensitivityadapted} to settings in which the models are required to satisfy a martingale coupling constraint. 
Our approach relies on a weaker version of the Implicit Function Theorem, which enables the construction of families of measures satisfying the prescribed constraints. 
We provide closed‑form expressions for these sensitivities, along with a characterization of the hedging strategies when the underlying process is real‑valued.
\end{abstract}


\vspace{3mm}

\section{Introduction}

\paragraph{Model Risk through Distributionally Robust Optimization.} 
Consider a criterion $g : \Pc(\X) \rightarrow \R$ for some space $\X$.  
In many contexts, we wish to evaluate $g$ at a measure $\mu$, which may represent a model chosen by the agent or determined through calibration methods. 
In either case, the agent inevitably faces uncertainty, either due to the multiplicity of plausible models or because calibration relies on imperfect information. 
This uncertainty, known as "Knightian uncertainty", introduced in \citeauthor{knight1921risk} \cite{knight1921risk}, has been widely studied in economics and decision theory.
In mathematical finance specifically, this model uncertainty is often referred to as model risk, and has notably been studied in the context of derivative pricing and robust hedging under uncertain volatility by \citeauthor{avellaneda1995pricing} \cite{avellaneda1995pricing}, \citeauthor{dupire1997pricing} \cite{dupire1997pricing}, and \citeauthor{hobson1998robust} \cite{hobson1998robust}.
A classical example arises in stochastic optimization, where the criterion takes the form $g(\mu) = \inf_{a \in A} \int_\X f(x, a)\mu(\mathrm{d}x)$.
In this setting, $\mu$ models the distribution of uncertainty in the environment, and the agent seeks an action $a \in A$ that minimizes expected cost.
Recently, Distributionally Robust Optimization (DRO) emerged as a systematic approach to deal with model uncertainty for stochastic optimization problems. 
In DRO, the agent competes against an adversary who perturbs the reference distribution within a prescribed deviation set $D \subseteq \Pc(\X)$.
This leads to the min-max formulation
\begin{equation}\label{ch03eqdef:droproblem}
\inf_{a \in A} \sup_{\mu' \in D} \int_\X f(x, a)\mu'(\mathrm{d}x).
\end{equation}
The choice of $D$ is central, as it encodes the extent and structure of model uncertainty. For comprehensive surveys of DRO and its applications, we refer to \citeauthor{rahimian_distributionally_2019} \cite{rahimian_distributionally_2019}, \citeauthor{kuhn2025distributionally} \cite{kuhn2025distributionally}, \citeauthor{lin2022distributionally} \cite{lin2022distributionally} for extensive reviews of DRO.
Among the many possible constructions of $D$, sets based on an optimal transport criterion have been a subject of interest and have been initiated in \citeauthor{blanchet2019quantifying} \cite{blanchet2019quantifying}.
In this paper we focus on the $p$-Wasserstein distance and its $p$-adapted variant, following the framework developed in \citeauthor{Touzisauldubois2024ordermartingalemodelrisk} \cite{Touzisauldubois2024ordermartingalemodelrisk}.
Two main approaches have been developed for analyzing DRO under Wasserstein-type uncertainty sets:
The DRO problem \eqref{ch03eqdef:droproblem} is reformulated as a finite-dimensional optimization problem. 
This approach has been applied in \citeauthor{zhang2025short} \cite{zhang2025short}, \citeauthor{ji2021data} \cite{ji2021data}, \citeauthor{blanchet2019quantifying} \cite{blanchet2019quantifying}, \citeauthor{mohajerin_esfahani_data-driven_2018} \cite{mohajerin_esfahani_data-driven_2018}, and for causal variants in \citeauthor{han2022distributionally} \cite{han2022distributionally}, \citeauthor{jiang2024duality} \cite{jiang2024duality}.
The other approach that we will use here is the sensitivity analysis.
Here, the deviation set $D$ is defined as a Wasserstein ball of radius $r$, and the behavior of~\eqref{ch03eqdef:droproblem} is studied asymptotically as $r \to 0$. This approach yields tractable expansions and insights into robustness. 
It has been developed in \citeauthor{bartl2021sensitivity} \cite{bartl2021sensitivity}, \citeauthor{bartlsensitivityadapted} \cite{bartlsensitivityadapted}, \citeauthor{bartl2023sensitivity} \cite{bartl2023sensitivity}, and \citeauthor{jiang2024sensitivity} \cite{jiang2024sensitivity}, and is also the perspective adopted in this paper.

\paragraph{Model Risk Hedging through Constrained DRO.} 
It is natural to study distributionally robust optimization (DRO) problems with constraints for two main reasons. 
First, from a practical perspective, it is unrealistic to discard all information used during model calibration, thus adding constraints to the family of deviations considered.
Furthermore, as we discussed in \citeauthor{Touzisauldubois2024ordermartingalemodelrisk} \cite{Touzisauldubois2024ordermartingalemodelrisk}, using available market instruments for hedging reduces risk exposure and thus helps manage market risk. 
By standard dualization of the corresponding hedging strategies, the robust evaluation reduces to considering neighboring models that simultaneously satisfy martingale and/or marginal constraints. 
The latter represents calibration to current market information provided by the implied volatility surface. 
Following the literature on martingale optimal transport \citeauthor{beiglbock2016problem} \cite{beiglbock2016problem}, \citeauthor{acciaio2016model} \cite{acciaio2016model}, we consider the following two cases:
\begin{itemize}
\item All returns generated by buy-and-hold strategies $\mathfrak{H}_{\rm M}:=\big\{h(X_1)\cdot(X_2-X_1):h\in C^0_b(\R^d)\big\}$;
\item All returns generated by a vanilla payoff with maturity $T\in\{1,2\}$ define the set of zero-cost strategies $\mathfrak{H}_{\rm m_T}:=\big\{f(X_T)-\mu_T(f):f\in C^0_b(\R^d)\big\}$, where $\mu_T=\mu\circ X_T^{-1}$ denotes the $T$th marginal of $\mu$. 
\end{itemize}
Of course, other instruments exist. 
For instance, one could also consider a $\rm{VIX}$ future or a forward $\rm{log}$ contract, as studied in \citeauthor{de2015linking} \cite{de2015linking}. 
However, we do not consider them in this article. 
Given a set $\mathfrak{H}$ of hedging instruments, and recalling that $\W_p^{\rm ad}$ denotes the adapted Wasserstein distance on $\Pc(\R^d \times \R^d)$ (defined in Section~\ref{sec:notations}), we define the upper and lower model distributionally robust hedging problems by
$$
\overline{G}^{\mathfrak{H}}(r)
:=
\inf_{\mathfrak{h} \in\mathfrak{H}} \sup_{ \W_p^{\rm ad}(\mu, \mu') \le r} \big\{ g ( \mu' )+\mu'(\mathfrak{h}) \big\}
~~\mbox{and}~~
\underline{G}^{\mathfrak{H}}(r)
:=
\sup_{\W_p^{\rm ad} (\mu, \mu') \le r} 
\inf_{\mathfrak{h}\in\mathfrak{H}} \big\{ g ( \mu' ) +\mu'(\mathfrak{h})\big\}
,$$
where $\mathfrak{H}$ is either $\mathfrak{H}_{\rm M}$, $\mathfrak{H}_{\rm m}:= \mathfrak{H}_{\rm m_1} + \mathfrak{H}_{\rm m_2}$, or $ \mathfrak{H}_{\rm M} + \mathfrak{H}_{\rm m}$.
The second problem is a constrained DRO. 
Specifically, when the instrument set is limited to buy-and-hold strategies, the problem falls under the scope of \citeauthor{Touzisauldubois2024ordermartingalemodelrisk} \cite{Touzisauldubois2024ordermartingalemodelrisk} and also \citeauthor{jiang2024sensitivity} \cite{jiang2024sensitivity}.
Note that \citeauthor{lam2018sensitivity} \cite{lam2018sensitivity} considered the case where the set of instruments consists of returns generated by vanilla payoffs, and the distance is replaced by the Kullback divergence. 
In \citeauthor{Touzisauldubois2024ordermartingalemodelrisk} \cite{Touzisauldubois2024ordermartingalemodelrisk} we left open the case where vanilla payoffs with maturity $2$ are available for a partial hedging purpose. 
This will be the object of this article.

\noindent \paragraph{Our contribution: A new method to study constrained DRO.}

It turns out that, in the $2$-period setting, studying the adapted Wasserstein DRO problem where deviations are martingale couplings is involved.
We will show that for hedging instruments consisting of semi-static and dynamic hedging strategies, 
$
\mathfrak{H}$
with $\E^\mu[\mathfrak{h}] = 0 $ for all $ \mathfrak{h} \in \mathfrak{H}$,
the derivatives at zero satisfy  
$$
\left.\overline{G}^{\mathfrak{H}}\right.'(0)
=
\left.\underline{G}^{\mathfrak{H}}\right.'(0)
=
\inf_{\mathfrak{h}\in\mathfrak{H}}
\big\|
\partial_x^{\rm ad} (g + \mathfrak{h})
\big\|_{\mathbb{L}^{p'}(\mu)},
$$
where $\partial_x^{\rm ad}$ is the differential operator associated with the adapted Wasserstein distance.
We propose a new approach for constrained DRO, based on the Implicit Function Theorem.
This method enables us to handle a broader class of constrained DRO problems than those considered in \citeauthor{jiang2024sensitivity} \cite{jiang2024sensitivity} and \citeauthor{Touzisauldubois2024ordermartingalemodelrisk} \cite{Touzisauldubois2024ordermartingalemodelrisk}.
Due to the already long and technical nature of this paper, the study of sensitivity for DRO under a finite number of constraints or just law constraints is currently a work in progress (\citeauthor{CrowellSaulduboisTouziLawConstraints} \cite{CrowellSaulduboisTouziLawConstraints} and  \citeauthor{CrowellSaulduboisTouziFiniteConstraints} \cite{CrowellSaulduboisTouziFiniteConstraints}).

\noindent The key idea behind this new method relies on finding a good approximation of
$\text{argmax}_{\mu', \W_p^{\rm ad} (\mu, \mu') \le r} \inf_{\mathfrak{h}\in\mathfrak{H}} \{g ( \mu' )+\mu'(\mathfrak{h}) \}$.
The natural candidate considered in \citeauthor{bartl2021sensitivity} \cite{bartl2021sensitivity}, \citeauthor{bartlsensitivityadapted} \cite{bartlsensitivityadapted} and in \citeauthor{Touzisauldubois2024ordermartingalemodelrisk} \cite{Touzisauldubois2024ordermartingalemodelrisk} has the form $\mu^{r T^{\mathfrak{h}}} := \mu\circ (X + rT^{\mathfrak{h}})^{-1}$ for some map $T^{\mathfrak{h}}$ defined through the minimizer $\text{argmin}_{\mathfrak{h}\in\mathfrak{H}} \Vert \partial_x^{\rm ad} ( \delta_ m g + \mathfrak{h} ) \Vert_{\mathbb{L}^{p'} (\mu) } $, where $\delta_ m g$ is the linear functional derivative of $g$ (see definition \ref{def:linearfunctionderiv}).
However, in contrast with the situation of \citeauthor{Touzisauldubois2024ordermartingalemodelrisk} \cite{Touzisauldubois2024ordermartingalemodelrisk}, this natural candidate does not generally satisfy the constraint associated with $\mathfrak{H}$.
To overcome this difficulty, we use a variation of the Implicit Function Theorem to find sharp estimates of
$$
\W_p^{\rm ad} ( \mu^{r T^{\mathfrak{h}}}, \mathfrak{H}^{\perp}) = 
\inf_{\substack{\nu \in \Pc(\X) \\ \text{for all $ \mathfrak{h} \in \mathfrak{H}, $ } \E^\nu[ \mathfrak{h}] = 0} }
\W_p^{\rm ad} (\nu,  \mu^{r T^{\mathfrak{h}}})
.$$
We finally illustrate our results in Section \ref{sec:num_illustration} and represent the first-order hedging strategies and the sensitivities for both Bachelier and Black-Scholes models. 
We find that adding the second marginal to the constraints significantly decreases the sensitivity. 
It is also worth noting that buy-and-hold hedging strategies are significantly impacted by the possibility of adding vanilla payoffs to the hedging strategy.

\section{Notations and Definitions}\label{sec:notations}

Throughout this paper, let $ p > 1$ be a real number and $p' = \frac{p}{p-1}$ denotes the conjugate exponent of $p$.
We denote $S:= \R^d$ and $\X:= S \times S$ for some integer $d\ge 1$, both endowed with the corresponding canonical Euclidean structure and the associated norm defined by $x\longmapsto\vert x \vert:=\sqrt{x \cdot x}$.
For $k \in \N$, let $ C^k_b(\X, \R)$ denote the space of functions from $\X$ to $ \R$ that are $k$-times continuously differentiable, with all derivatives up to order $k$ being bounded.
We write an element of $\X$ as $ x:= (x_1, x_2) \in \X$, and define
\begin{equation}\label{eqdef:n and nad}
 \mathbf{N} (x)  := \frac{\nabla(\vert \cdot \vert^{p'}) (x) }{p'}
 = \frac{x}{ \vert x \vert^{2 - p'}},
 \,\,
 \Nad (x) :=
 \begin{bmatrix}
 \mathbf{N} (x_1)
\\[3pt]
 \mathbf{N} (x_2)
\end{bmatrix}
.
\end{equation}
We intentionally abuse the notation $\mathbf{N} (x_i)$ and $\mathbf{N} (x)$ as the space on which $ \mathbf{N} $ is defined is implied by the context of the variable.
Note that, for $ d = 1$, $ \mathbf{N} (x) = \text{sgn} (x) \vert x \vert^{p' - 1}$.
\begin{Definition}\label{def:coercivity}
Let $(E, \vert \cdot \vert )$ be a normed vector space. A function $ f: E \rightarrow \R $ is said to be coercive if it satisfies $ f(x)  \xrightarrow[\vert x \vert \rightarrow +\infty]{} + \infty$.
\end{Definition}

\noindent Let $\Pc (\X) $ be the collection of all probability measures $\mu$ on $\X$.
We denote by $ \Pc_p (\X) $, the subset of those with finite $p$-th moment:
$$
\Pc_p (\X):= \{\mu\in\Pc (\X):~\E^\mu[|X|^p]<\infty \}.
$$
For a subset $ \Ec \subset C^0_b(\X, \R)$ (bounded continuous functions), we define
\begin{equation}\label{eqdef:orthogonalityduality}
\Ec^{\perp}:= \{ \mu' \in \Pc(\X) \,: \,\int f \mathrm{d}\mu' = 0 \, \, \text{for all } f \in \Ec \}.
\end{equation}
We define the projection maps $(X,X')$ in $\X\times \X$ defined by $X(x,x')=x$ and $X'(x,x')=x'$ for all $x,x'\in \X$.
For $\mu,\mu'\in\Pc (\X ) $, we define the set of all couplings
$$
\Pi(\mu,\mu')
:=
\{ \pi\in\Pc (\X \times \X):\pi\circ X^{-1}=\mu~\mbox{and}~\pi\circ {X'}^{-1}=\mu' \}.
$$
In the dynamic setting, we extract the time components using the projection maps $X_i$, $X_i'$, defined on $\X \times \X$:
$$
X_i (x, x') = x_i
~\mbox{and}~
X'_i (x,x') = x'_i,
~i=1, 2,
~\mbox{for all}~
(x, x')\in \X\times\X.
$$
We define $L^p(\mu):= \{\text{$f$ measurable, satisfying} \int_{\X} \vert f \vert^p \mathrm{d}\mu < \infty  \}$ and $\mathbb{L}^p(\mu)$ is then defined as the quotient of $L^p(\mu)$ by equivalence relation $\mu$-almost everywhere equality.
Note that $ \mathbb{L}^p(\mu)$ is a Banach space while $L^p(\mu)$ is not.
We define $\mathbb{L}^p_0(\mu)$ as the subset of $0$-mean functions in $\mathbb{L}^p(\mu)$:
\begin{equation}\label{eqdef:L0}
\mathbb{L}^{p}_0(\mu) := \{ f \in \mathbb{L}^p (\mu) : \, \E^{\mu}[f(X)] = 0 \}.
\end{equation}
For an element $T : \X \rightarrow \X$ of $\mathbb{L}^{p}(\mu)$, define the pushforward measure $\mu^{ T} \in \Pc(\X)$ as
\begin{equation}\label{eqdef:pushforward}
    \mu^{T} := \mu \circ ( X + T )^{-1}.
\end{equation}

\begin{Definition}
A probability measure $ \mathbb{P} \in \mathcal{P}(\X \times \X) $ is causal if the filtration generated by $X$, $\mathbb{F}^{X}$ is compatible with the filtration generated by $X'$ : $\mathbb{F}^{X'}$, in the sense that for all bounded Borel-measurable $ f: \X \rightarrow \R $ and $ g: S \rightarrow \R $,
$$
\E^{\P} \big[
f (X_1, X_2) g (X'_1) \vert X_1 \big]
=
\E^{\P} \big[
f (X_1, X_2) \vert X_1 \big]
\,
\E^{\P} \big[
g (X'_1) \vert X_1 \big]
.
$$

\end{Definition}

We introduce the set of {\it bi-causal} couplings
$$
\Pi^{\text{bc}}(\mu,\mu')
:=
\{\pi\in \Pi (\mu, \mu') \,\, \text{such that} \,\, \pi \,\, \text{and} \,\,\pi \circ (X', X)^{-1} \,\, \text{are causal} \},
$$
together with the corresponding {\it adapted Wasserstein} distance
\begin{eqnarray*}
\W^{{\rm ad}}_p(\mu,\mu')
:=
\inf_{\pi\in\Pi^{\text{bc}}(\mu,\mu')}\E^\pi\Big[\big|X-X'\big|^p\Big]^{\frac1p},
&\mbox{for all}&
\mu,\mu'\in\Pc_p(\X).
\end{eqnarray*}
This causal/bicausal and adapted Wasserstein framework was introduced and developed in \citeauthor{backhoff2017causal} \cite{backhoff2017causal}, \citeauthor{lassalle2018causal} \cite{lassalle2018causal}, and \citeauthor{backhoff_adapted_2020} \cite{backhoff_adapted_2020}; see also \citeauthor{blanchet2024empirical} \cite{blanchet2024empirical} for a recent use of this distance under a martingale constraint.
We denote the corresponding ball of radius $r$ by
$$
\Bad (\mu,r)
:=
\{\mu'\in\Pc_p(\X):\W^{{\rm ad}}_p(\mu,\mu')\le r \}
.$$
Here and throughout the article, the continuity of a map defined on $\Pc_p$ refers to the adapted Wasserstein distance $\W_p^{\rm ad}$.
We introduce the following notation.
For $ u \in C^{1} (\X, \R) $ we set for $i = 1,2$, $ \E^{\mu}_i [ u ]:= \E^{\mu}[ u(X) \vert X_i ] $, and define the corresponding gradient
$$
\partial_x^{\rm ad} u:= \begin{pmatrix}
\E^{\mu}_1 [ \partial_{x_1} u ]
\\
\partial_{x_2} u
\end{pmatrix}
.
$$

 \begin{Definition}\label{def:adaptedfunctions}
     We say that $T : \X \rightarrow \X$ is adapted if and only if, for all $x \in \X$, $ T (x) = (T_1(x_1), T_2(x_1, x_2) )$
 \end{Definition}

 \noindent In the following, we consider $\mathbb{L}^p_{\rm ad} ( \mu )$
 \begin{equation}\label{eqdef:adapted Lp}
\mathbb{L}^p_{\rm ad} ( \mu ) := \{ T \in \mathbb{L}^p ( \mu ) \,\, \text{such that }  T \, \text{is adapted} \}
 \end{equation}
 endowed with the norm $ \Vert T \Vert^p_{\mathbb{L}^p_{ \rm ad}  ( \mu )} := \Vert T_1 \Vert_{\mathbb{L}^p ( \mu_1)}^p + \Vert T_2 \Vert_{\mathbb{L}^p ( \mu)}^p  $.
Finally, for a subset $U \subset \Pc_p(\X)$, define the distance between $\mu$ and $U$ as
\begin{equation}\label{eqdef:distance to a part}
\W_p^{\rm ad} ( \mu, U ) := \inf_{ u \in U}\W_p^{\rm ad} ( \mu, u )
.\end{equation}
Throughout this paper, we consider a function $g:\Pc_p(\X)\longrightarrow\R$ with appropriate smoothness in the following sense.
A function is said to have $p$-polynomial growth if it is bounded uniformly by $C(1+\vert x\vert ^p)$ for some constant $C$.
\begin{Definition}\label{def:linearfunctionderiv}
We say that $g$ has a linear functional derivative if there exists a continuous function $\delta_mg:\Pc_p(\X) \times \X \longrightarrow\R$, with $p$-polynomial growth in $x$, locally uniformly in $\mu$, such that for all $\mu,\mu'\in\Pc_p(\X)$, and denoting $\bar\mu^\lambda:=\mu+\lambda(\mu'-\mu)$, we have
$$
\frac{g (\bar\mu^\lambda )-g(\mu)}{\lambda}
\longrightarrow
\langle\delta_mg(\mu,x),\mu'-\mu\rangle
:=
\int_\X \delta_mg(\mu,x)(\mu'-\mu)(dx),
~\mbox{as}~\lambda\searrow 0.
$$
\end{Definition}

Clearly, the linear functional derivative is defined up to a constant, which will be irrelevant throughout this paper.
In the linear case where $ g (\mu) = \int_{\X} f \mathrm{d} \mu$, for some continuous map $f$ with $p$-polynomial growth, the linear functional derivative is the constant map (in $\mu$) $ \delta_m g (\mu, x) = f(x) $ for all $\mu \in \Pc_p \big(\X \big) $, $ x \in \X $.
We recall that, up to technical conditions, the Lions' derivative coincides with the Wasserstein gradient. It is then given by $\partial_x\delta_m g$, see \citeauthor{CarmonaDelarue} \cite{CarmonaDelarue}.
Moreover, this definition is equivalent to the existence of such a continuous function $\delta_mg$ satisfying:
\begin{equation*}
g(\mu')-g(\mu)
=
\int_0^1 \int_\X \delta_mg(\bar\mu^\lambda,x)(\mu'-\mu)(dx)d\lambda,
\,\,\, \mbox{for all} \,\,\,
\mu,\mu'\in\Pc_p (\X).
\end{equation*}

\section{Main Results}

We study the sensitivity of distributionally robust optimization (DRO) with respect to the adapted Wasserstein distance under a martingale coupling constraint, \textit{i.e.}, when the competing models are required to share the same two marginals as $\mu$ and to remain martingale measures.
We then adapt the result to American options via optimal stopping problems.

\begin{Assumption} \label{ass:regul and estim on g}
The mapping $g$ has a linear functional derivative such that $\delta_mg$ is $C^1$ in $x$, and $\partial_x\delta_mg$ is jointly continuous, with $(p-1)$-polynomial growth in the $x$-variable locally in the $\mu$-variable.
\end{Assumption}

\begin{Remark}\label{rem:existence normalisation and distance}
{\rm
Throughout this article, whenever we normalize a map $T$ of the form $T = \frac{1}{c} \Nad (\Lambda)$, for some $\Lambda \in \mathbb{L}^{p'} (\mu)$, so that $\Vert T \Vert_{\mathbb{L}^p (\mu)} = 1$, the normalizing constant $c$ exists and is unique whenever $\Lambda \neq 0$ in $\mathbb{L}^{p'} (\mu)$.
The degenerate case $\Lambda = 0$ corresponds to the infimum defining the associated sensitivity already being attained at $0$, so that the differentiability of the corresponding upper and lower DRO maps is then immediate; see Remark $6.3$ of \citeauthor{Touzisauldubois2024ordermartingalemodelrisk} \cite{Touzisauldubois2024ordermartingalemodelrisk}.
Here and in the rest of the article, we implicitly assume that $\Lambda \neq 0$, so that such a constant $c$ always exists.
}
\end{Remark}

\subsection{Semi-Static Hedging Strategy}\label{subsec:semi static hedging}

For $\mu \in \Pc_p (\X)$, the corresponding marginals are denoted $ \mu_i = \mu \circ X_i^{-1} $ for $ i=1, 2$.
We denote by $F_2 (x) := \mu_2 ( (-\infty, x] )$ the cumulative distribution function of $\mu_2$ and by $F_2^{-1} (u) := \inf \{ x \in \R : F_2 (x) \geq u \}$, $u \in (0, 1)$, its generalized inverse.
We set:
\begin{equation}\label{eqdef:marginals}
\Pi_i (\mu_i) := \{ \mu' \in \Pc(\X) \, : \, \mu'\circ X_i^{-1} = \mu_i \, \}
\text{ for $i =1, 2$ }
,\end{equation}
along with $\Pi(\mu_1, \mu_2) :=  \Pi_1 (\mu_1) \cap \Pi_2 (\mu_2)$, and, for $x = (x_1, x_2) \in \R \times \R$,
\begin{equation}\label{eqdef:oplusnotation}
f_1 \oplus f_2 (x) = f_1(x_1) + f_2(x_2).
\end{equation}

Now, we allow semi-static hedging strategies, \textit{i.e.}, one can hold a vanilla payoff $f_1 (S_1)$ and $f_2 (S_2) $ and dynamically hedge one's position (see \cite{henry2017model} for more details).
We define ${\rm M} $ as the set of all probability measures on $\X$ that are martingales, 
\begin{eqnarray*}
\rm M
&:=&
\{\mu\in\Pc (\X):~\E^\mu[X_2|X_1]=X_1,~\mu-\mbox{\textit{a.s.}} \},
\end{eqnarray*}
and recall that for $\mu\in \Pc_p (\X)$, we have
\begin{eqnarray*}
\mu\in \rm M
&\mbox{if and only if}&
\E^\mu\big[h^\otimes \big]=0,
~\mbox{for all}~
h\in \mathbb{L}^{\infty}(\mu_1),
\end{eqnarray*}
where we used the notation 
\begin{equation}\label{eqdef:otimes notation}
h^\otimes(x)
:=h(x_1)\cdot(x_2-x_1),
\text{ for all }
x=(x_1,x_2)\in \X,
\end{equation}
and we also use the notation $f_1 \oplus f_2$ defined by \eqref{eqdef:oplusnotation}.
Define $\Pi^{\rm M} (\mu_1, \mu_2) := {\rm M} \cap \Pi (\mu_1, \mu_2)$ (see Definitions \eqref{eqdef:marginals}), and the corresponding constrained ball $\BadMm (\mu, r) := \Bad (\mu, r) \cap \Pi^{\rm M} (\mu_1, \mu_2)$.
We will extensively use the identification between $\mathbb{L}^p (\mu_i)$ and its canonical injection in $\mathbb{L}^p(\mu)$: for $f \in \mathbb{L}^p(\mu_i)$ (where $f$ is a representative of the equivalence class), we identify $f$ with the map $x \in \X \mapsto f(x_i)$, which belongs to $\mathbb{L}^p(\mu)$ since $\mu\circ X_i^{-1} = \mu_i$.
Due to the complexity of the DRO problem under Wasserstein and martingale constraints, here we consider the adapted Wasserstein metric. 
We define the following two functions
\begin{equation*}
\upGmartcoupad (r) :=   \inf_{ \substack{ f_1, f_2, h \in C^1_b }}
															 \sup_{ \mu' \in \Bad (\mu, r)}    g(\mu') + \int  ( f_1  \oplus  f_2  +  h^{\otimes}) \mathrm{d} (\mu'  - \mu) 
\end{equation*}
and
$
\lowGmartcoupad (r) :=      \sup_{ \mu' \in \BadMm (\mu, r)}
															 g(\mu') 
$.

\begin{AssumptionC} \label{ass:mart coupling}
Let $ \mu \in \Pc_p (\R^2)$ be a martingale measure.
\begin{enumerate}[label=\textnormal{(\roman*)}]
\item \label{cond:disintegration}$\mu$ admits the disintegration $ \mu(\mathrm{d}x):= \kappa(x_1, x_2) \mathrm{d}x_2 \mu_1 (\mathrm{d} x_1) $, with $\kappa \in \mathbb{L}^{\infty} (\R^2)$ and $ \kappa (X) \geq c >0 $ $\mu$-almost surely, for some constant $c >0$.
\item\label{cond:density second marginal} $\mu_2 (\mathrm{d} x _2) = q_2(x_2) \mathds{1}_{ \{ x_2 \in I_2 \}} \mathrm{d}x_2$ for some bounded interval $I_2:= [a, b ]$, and some strictly positive continuous $q_2$.
\item\label{cond:info discrepancy} The following informational discrepancy holds: $ \mathbb{L}^1(\mu_1) \cap \mathbb{L}^1(\mu_2) = \R$. 
\end{enumerate}
\end{AssumptionC}

\begin{Remark}
    {\rm 
 \textbf{In Assumption \hyperref[ass:mart coupling]{$\textbf{\rm C}_{\W_p^{\rm ad}}$} \ref{cond:info discrepancy}, the spaces $ \mathbb{L}^1(\mu_1)$ and $ \mathbb{L}^1(\mu_2)$
 are understood through their canonical injections into $\mathbb{L}^1(\mu)$.}
 Moreover $\R$ is identified with the set of $\mu$-almost-everywhere constant functions. We discuss this assumption in 
the Appendix \ref{subsec:discussion}.}
\end{Remark}

\noindent We begin by estimating the adapted Wasserstein distance between $\mu^{r \Theta}$ and $ \Pi^{\rm M} (\mu_1, \mu_2) $. 
Define the following operator
\begin{equation}\label{eqdef:operatormartcouplings}
   \opedistmartcoup: \Theta \in \mathbb{L}^p_{\rm ad} (\mu) \mapsto
   \big( \E^{\mu}_1 [ \Theta_1 ], \E^{\mu}_2[\Theta_2], \E^{\mu}_1 [\Theta_2 - \Theta_1 ] \big) \in \prod_{i=1}^2 \mathbb{L}^p (\mu_i)  \times  \mathbb{L}^p_0 (\mu_1 )
.\end{equation}

\begin{Proposition}\label{prop:proj mart coupling set}
Let $\mu$ satisfy Assumption \hyperref[ass:mart coupling]{$\textbf{\rm C}_{\mathbb{W}^{\rm ad}_p}$}. Let $ \Theta_2: \R^2 \rightarrow \R$ be a compactly supported $C^2$ function. 
Set $ \Theta := ( 0, \Theta_2)$, then, there exists $ C > 0 $ such that
\begin{equation}
\limsup_{r \rightarrow 0} \frac{1}{r}
\W^{\rm ad}_p \big( \mu^{r \Theta},  \Pi^{\rm M} (\mu_1, \mu_2) \big)
\leq 
C
\Vert \opedistmartcoup  (\Theta ) \Vert,
\end{equation}
where $\opedistmartcoup$ is defined by \eqref{eqdef:operatormartcouplings}.
\end{Proposition}

\begin{Proposition}\label{prop:sensi 1 dim mart coupling} 
Let $g$ satisfy Assumption \ref{ass:regul and estim on g} and $\mu$ satisfy \hyperref[ass:mart coupling]{$\textbf{\rm C}_{\W_p^{\rm ad}}$}. 
Then both maps $\lowGmartcoupad$ and $\upGmartcoupad$ are differentiable at $0$, and 
\begin{equation}\label{eq:deriv0 martingale coupling}
 \left. \upGmartcoupad \right.'   (0) 
=
 \left. \lowGmartcoupad \right. '  (0) 
=
\inf_{ \substack{ (f, h) \in ( \Pi_{i =1}^2 \mathbb{L}^{p'} (\mu_i) ) \times  \mathbb{L}^{p'}_0 (\mu_1)}}
\left. U^{ \rm M, m}_{ \rm ad} (f , h)  \right.^{1/p'}
,\end{equation}
where $
U^{ \rm M, m}_{ \rm ad} (f, h)
:=
 \Vert
\partial_x^{\rm ad} \delta_m g + f + h J
\Vert^{p'}
_{p'}$ and $J := (-1, 1)^{\intercal}$, and where $\mathbb{L}^{p'}_0(\mu_1)$ is defined by \eqref{eqdef:L0}.
Restricting $h$ to $\mathbb{L}^{p'}_0 (\mu_1)$ is no loss of generality: for $h \in \mathbb{L}^{p'} (\mu_1)$, the constant part $\E^{\mu}[h] J = (- \E^{\mu}[h], \E^{\mu}[h])^{\intercal}$ can be absorbed into $f = (f_1, f_2)$, since constants belong to both $\mathbb{L}^{p'} (\mu_1)$ and $\mathbb{L}^{p'} (\mu_2)$.
Then, $ U^{ \rm M, m}_{ \rm ad}: \mathbb{L}^{p'} (\mu_1) \times \mathbb{L}^{p'} (\mu_2) \times \mathbb{L}^{p'}_0 (\mu_1) \rightarrow \R $ is strictly convex, continuous, and coercive (see Definition \ref{def:coercivity}). Hence, the optimization problem \eqref{eq:deriv0 martingale coupling} admits a unique minimizer $ (f^{ {\rm M, m}}, h_{\rm M, m}) \in  \Prod_{i=1}^2 \mathbb{L}^{p'} (\mu_i) \times \mathbb{L}^{p'}_0 (\mu_1) $, characterized by the first-order equations:
\begin{equation}\label{eqdef:FOC Mart couplings}
\opedistmartcoup (T^{\rm M, m}) = 0 
\,\, \text{where} \,\, 
 T^{\rm M, m} := \frac{1}{c}\Nad (\partial_{x}^{\rm ad} \delta_m g + f^{{\rm M, m}} + h_{\rm {M, m}}J), 
 \end{equation}
where $\opedistmartcoup$ is defined by \eqref{eqdef:operatormartcouplings} and $c$ is uniquely defined by $\Vert T^{\rm M, m}\Vert_{\mathbb{L}^p(\mu)} = 1$. \footnote{see Remark \ref{rem:existence normalisation and distance} for the existence of $c$}

\noindent Furthermore, when $p = 2$, the optimal strategies are given by 
\begin{equation*}
\begin{split}
f^{\rm M, m} _{1} = h_{\rm M, m} - \E^{\mu}_1 [\partial_{x_1} \delta_m g ] \,\, , \,\, 
f^{\rm M, m}_{ 2} = - \E^{\mu}_2 [ \partial_{x_2} \delta_m g + h_{\rm M, m} ]
,\end{split}
\end{equation*}
where $ h_{ \rm M, m} $ is the unique (up to a constant) solution of the following Fredholm equation:
$$
h_{ \rm M, m} - \E^{\mu}_1 [ \E^{\mu}_2 [ h_{ \rm M, m} ] ] = \E^{\mu}_1 [ \E^{\mu}_2 [ \partial_{x_2} \delta_m g ] ] - \E^{\mu}_1 [ \partial_{x_2} \delta_m g ] 
.$$
\end{Proposition}

\medskip

\begin{Remark}
{\rm
\noindent $\bullet$ Assumption \hyperref[ass:mart coupling]{$\textbf{\rm C}_{\W_p^{\rm ad}}$} \ref{cond:info discrepancy} is here to ensure that $\mu$ satisfies the following property that for $f \in \mathbb{L}^1 (\mu_1)$ and $g \in \mathbb{L}^1 (\mu_2) $ we have 
$$
f(X_1) = g(X_2) \,\, \mu-\textit{ a.s.} \implies \exists c \in \R \,\, \text{such that } f = c \,\, \mu_1-\textit{a.s} 
\,\, \text{and} \,\, 
 g = c \,\, \mu_2-\textit{a.s} 
.$$
As it will be discussed in Appendix \ref{subsec:discussion}, it turns out to be a condition on the support of the measure.

\noindent $\bullet$ Assumption \hyperref[ass:mart coupling]{$\textbf{\rm C}_{\W_p^{\rm ad}}$} \ref{cond:info discrepancy} is here to ensure the coercivity of the objective function of problem \eqref{eq:deriv0 martingale coupling}. 
In fact, as will be shown in the proof section, this assumption provides a sufficient condition for the contraction property of the operator $ \E^{\mu}_1 \circ \E^{\mu}_2: \mathbb{L}^\alpha_0 (\mu_1) \rightarrow \mathbb{L}^\alpha_0 (\mu_1)$ for all $1 \leq \alpha \leq \infty $. 
\textbf{However, we emphasize that this assumption is restrictive.} 
For instance, it does not include backward Monge couplings, which are proven to be dense in \citeauthor{nutz2024martingale} \cite{nutz2024martingale} (a backward Monge coupling has the form $\mu_2 \circ (T, \mathrm{Id})^{-1}$, \textit{i.e.}, it is the law of $(T(Y), Y)$ with $Y \sim \mu_2$).

\noindent $\bullet$ The $d$-dimensional case is difficult to tackle using this approach. 
In fact, our proof will use a linearization result for the one-dimensional Monge-Ampère equation between the measure $ \mu_2 $ and $ \mu\circ(X_2 + r \Theta + a(X_1)) ^{-1}$, for $ a \in \mathbb{L}^p (\mu_1) $ in a neighborhood of $0$ and $ r > 0 $ close to $0$. 
Such linearization turns out to be difficult in higher dimension.
The martingale optimal transport problem itself is considerably more involved in dimension $d \geq 2$, where the structure of optimal plans and the associated duality are studied in \citeauthor{de2018quasi} \cite{de2018quasi}, \citeauthor{de2018local} \cite{de2018local}, \citeauthor{de_march_multidimensional_2018} \cite{de_march_multidimensional_2018}, and \citeauthor{de2019irreducible} \cite{de2019irreducible}; see also \citeauthor{acciaio2020causal} \cite{acciaio2020causal} and \citeauthor{beiglbock2018denseness} \cite{beiglbock2018denseness} for the related causal/adapted transport theory in continuous time.
}
\end{Remark}

\subsection{Optimal Stopping Problem}

The result of Proposition \ref{prop:sensi 1 dim mart coupling} can be extended to optimal stopping problems, \textit{i.e.}, to the case where
\begin{equation}\label{eqdef:O.S problem}
g(\mu'):= \inf_{ \tau \in \text{ST}} \E^{\mu'} [ \ell_\tau ],
\end{equation}
where $\text{ST}$ is the set of all stopping times with respect to the canonical filtration.
\begin{Assumption}\label{ass:O.S}
\begin{enumerate}[label=\textnormal{(\roman*)}]
\item The map $\ell: \X \times \left\{ 1, 2 \right\} \rightarrow \R $ is adapted, with maps $\ell_1,\ell_2$ continuously differentiable and $(p-1)$-polynomially growing.
\item The Optimal Stopping Problem \eqref{eqdef:O.S problem} admits a unique solution $\hat\tau$.
\end{enumerate}
\end{Assumption}

\begin{Proposition}\label{prop:optimal stopping}
Let $g$ be defined by \eqref{eqdef:O.S problem} with $\ell$ satisfying Assumption \ref{ass:O.S}.
Let $d = 1$, and let $\mu $ satisfy Assumption \hyperref[ass:mart coupling]{$\textbf{\rm C}_{\W_p^{\rm ad}}$}.
Then both maps $\upGmartcoupad $ and $ \lowGmartcoupad$ are differentiable at $0$ and:
\begin{equation*}
\left. \lowGmartcoupad \right.' (0)
=
\left. \upGmartcoupad \right.' (0)
=
\inf_{ \substack{ h \in \mathbb{L}^{p'}_0(\mu_1) \\ f \in \mathbb{L}^{p'}(\mu_1)\times  \mathbb{L}^{p'}(\mu_2)}}
\| \partial_{x}^{\rm ad } \ell_{\hat\tau} + f + h J \|_{\mathbb{L}^{p'}(\mu)},
\end{equation*}
where $J$ is defined in Proposition \ref{prop:sensi 1 dim mart coupling}.
\end{Proposition}

\section{Numerical illustration}\label{sec:num_illustration}

In the same spirit as \citeauthor{Touzisauldubois2024ordermartingalemodelrisk} \cite{Touzisauldubois2024ordermartingalemodelrisk}, we will compare sensitivities for the American put under the $2$-adapted Wasserstein distance. 
Throughout this section, $p = 2$.
Besides the martingale coupling case studied in this article, we define the deviation without constraints and the deviation under the marginal constraints only (without the martingale constraint):
$$ G_{\rm ad} (r):= \sup_{ \mu' \in B_{\W^{\rm ad}_2} (\mu, r)} g(\mu') ,
\text{ and} \,\, G^{\rm m}_{\rm ad} (r):= \sup_{ \mu' \in B_{\W^{\rm ad}_2} (\mu, r) \cap \Pi (\mu_1, \mu_2)} g(\mu')
.$$
By \citeauthor{bartlsensitivityadapted} \cite{bartlsensitivityadapted}, $G_{\rm ad}$ is differentiable at $0$ and
$
G_{{\rm ad}}' (0) = \|\partial_x^{\rm ad} \delta_m g \|_{\mathbb{L}^2(\mu)}
$,
and, by Propositions \ref{prop:sensi 1 dim mart coupling} and \ref{prop:optimal stopping}, the sensitivity under the martingale coupling constraint is given by \eqref{eq:deriv0 martingale coupling} with $p = p' = 2$.
For the marginal constraints alone, we admit the sensitivity formula
$
\left. G^{\rm m}_{{\rm ad}} \right.' (0) = \inf_{f} \|\partial_x^{\rm ad} \delta_m g + f \|_{\mathbb{L}^2(\mu)}
$,
where the infimum is taken over $f \in \mathbb{L}^{2}(\mu_1) \times \mathbb{L}^{2}(\mu_2)$; this is the analogue of \eqref{eq:deriv0 martingale coupling} without the martingale term $h J$.
As in \citeauthor{Touzisauldubois2024ordermartingalemodelrisk} \cite{Touzisauldubois2024ordermartingalemodelrisk}, we will compare sensitivities and relative sensitivities, given by the sensitivity divided by the price. 
The two families of models will be defined through a pair $(Z_1,Z_2)\sim \Nc (0, 1)\otimes\Nc (0, 1)$:
\begin{itemize}
\item the Black-Scholes model $ \mu_{\text{BS}}^{\sigma}:= \mathcal{L} \big(e^{-\frac{\sigma^2}{2} + \sigma Z_1}, e^{ -\sigma^2+ \sigma (Z_1 + Z_2)} \big)$,
\item and the Bachelier model $ \mu^{\sigma} _{\text{Bach}}:= \mathcal{L} \big(\sigma Z_1, \sigma (Z_1 + Z_2) \big) $.
\end{itemize}
Since those families of models are parametric, we also compare the sensitivities with the corresponding {\it Vega} $:= \partial_\sigma g(\mu^\sigma) $.

In this section, we consider the example of an American put option with intrinsic values $ \ell_t (X) = (e^{- \rho t} K - X_t )^+ $, $t=1,2$, with $K = 1.3$ and $\rho = 0.05$.

\subsubsection{Comparison of Sensitivities}

\begin{figure}[H]
\centering
\begin{subfigure}[b]{0.49\linewidth}
\centering
\includegraphics[width=\linewidth]{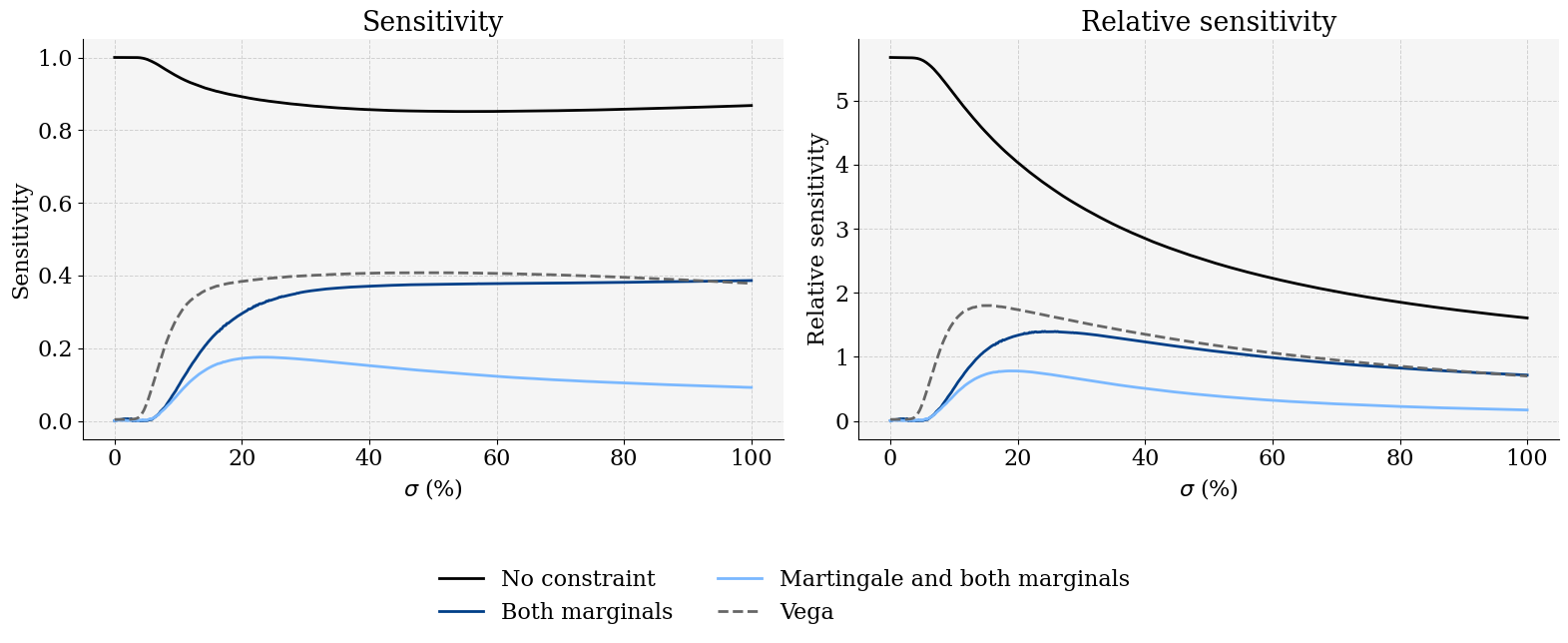}
\caption{\footnotesize Black-Scholes model.}
\label{fig:sensiBS}
\end{subfigure}
\hfill
\begin{subfigure}[b]{0.49\linewidth}
\centering
\includegraphics[width=\linewidth]{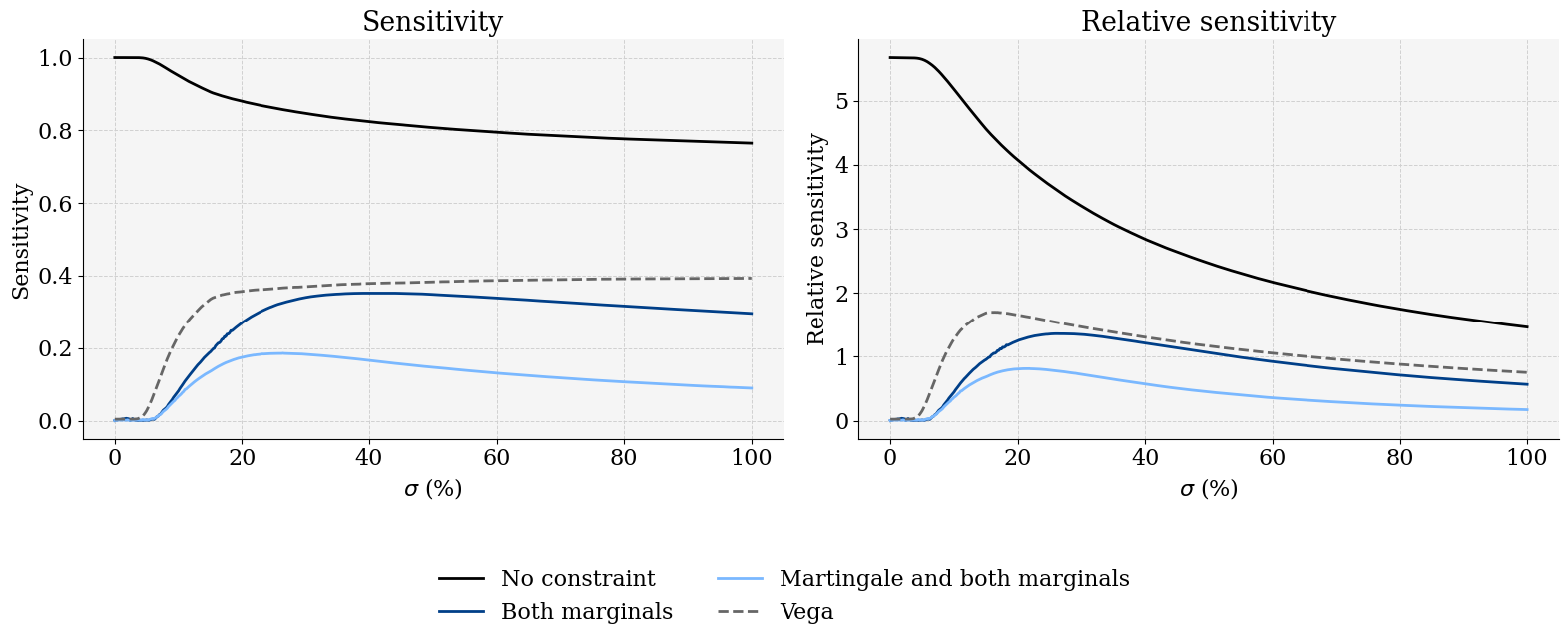}
\caption{\footnotesize Bachelier model.}
\label{fig:sensiBach}
\end{subfigure}
\caption{\it \footnotesize Sensitivities and \textit{Vega} for the buyer's price.}
\label{fig:sensi}
\end{figure}

We consider the sensitivity for the buyer's price $g(\mu') := \inf_{ \tau \in \text{ST}} \E^{\mu'} [ \ell_\tau ]$; the seller's price $g(\mu') := \sup_{ \tau \in \text{ST}} \E^{\mu'} [ \ell_\tau ]$ exhibits qualitatively similar behavior and is omitted.
Figure~\ref{fig:sensi} below compares, for the Black-Scholes and Bachelier models, the sensitivity without constraint, under the marginal constraints alone, and under the full martingale coupling constraint studied in this article, together with the \textit{Vega}.

We see a significant impact of adding the martingale coupling constraint on the sensitivity, as the discrepancy with the unconstrained and marginal-only cases is really noticeable, and a non-negligible difference with the \textit{Vega}: even though the \textit{Vega} evaluates model risk in a parametrized (martingale) family of measures, it does not account for the marginal constraints.

\subsubsection{Comparison of Hedging Strategies}

We can derive the hedging strategies from the computation of sensitivities.
We first consider hedging strategies for the
$$
\text{buyer's price :
$g(\mu') 
:=
\inf_{ \tau \in \text{ST}} \E^{\mu'} [ \ell_\tau ]
$
.}
$$

\begin{figure}[H]
\centering
\begin{subfigure}[b]{0.49\linewidth}
\centering
\includegraphics[width=\linewidth]{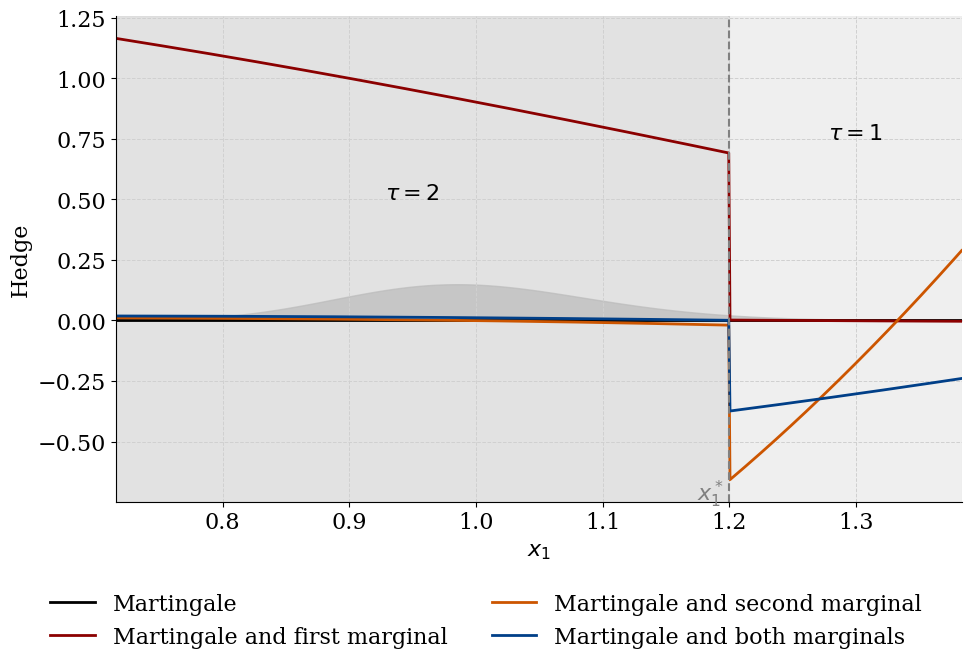}
\caption{\footnotesize $\sigma = 10\%$.}
\label{fig:hedgeamericanputBSsigma10}
\end{subfigure}
\hfill
\begin{subfigure}[b]{0.49\linewidth}
\centering
\includegraphics[width=\linewidth]{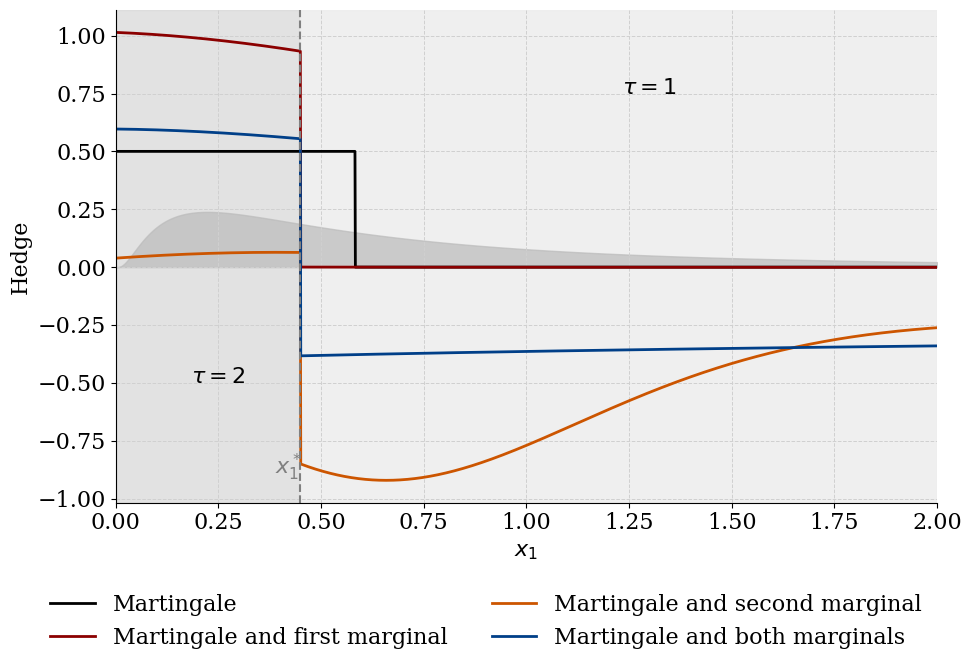}
\caption{\footnotesize $\sigma = 100\%$.}
\label{fig:hedgeamericanputBSsigma100}
\end{subfigure}
\caption{\it \footnotesize Hedging strategies in the \textbf{Black-Scholes} model.}
\label{fig:hedgeamericanputBS}
\end{figure}

\noindent We compare the case $\sigma = 10\%$ and a higher volatility $\sigma = 100\%$ for the Black-Scholes model (Figure~\ref{fig:hedgeamericanputBS}).

The jump visible in panel (a) reflects that the hedging strategy differs depending on whether the option holder enters the contract at $t = 1$, \textit{i.e.}, whether $\tau = 1$; at $\sigma = 10\%$ the probability of $\tau = 1$ is close to zero, while panel (b) shows that the hedging strategy becomes more involved at higher volatility.

\section{Proofs}

\subsection{A Sufficient Condition} \label{subsection:sufficientcondition}

We first provide a sufficient condition for differentiability at $0$ of a general constrained DRO problem under the adapted Wasserstein distance $\W_p^{\rm ad}$.
Let $\varphi : \Pc_p (\X) \rightarrow \R^k$, and let $\Econsgen \subset C^1_b (\X, \R) $ be a vector space.
Letting $\Econsgen^{\perp}$ be defined by \eqref{eqdef:orthogonalityduality}, consider the DRO problems
$$
\upGgenconsddist (r) := \inf_{ (\lambda, f ) \in \R^k \times \Econsgen} \sup_{ \mu' \in \Bad (\mu, r) } \Big\{ g(\mu') + \lambda \cdot \varphi(\mu') + \int f \mathrm{d}\mu' \Big\},
$$
and
$
\lowGgenconsddist (r) := \sup_{ \mu' \in \Bdisgen (\mu, r) } g(\mu')
$.
Here, $\Bdisgen (\mu, r):= \Bad (\mu, r) \cap \varphi^{-1} (\{ 0\}) \cap \Econsgen^{\perp}$.
We introduce the following set:
\begin{equation}\label{eqdef:closure image of gradient}
\text{$\clEconsdist$ is the $\mathbb{L}^{p'}(\mu)$-closure of the set $\partial_x^{\rm ad} (\Econsgen) = \{ \partial_x^{\rm ad} u \,\, : \,\, u \in \Econsgen \}$.}
\end{equation}
Finally, for $ \varphi, g$ which satisfy Assumption \ref{ass:regul and estim on g}, $ \lambda \in \R^k$ and $ u $ in $\clEconsdist$, define the following mapping
\begin{equation}\label{eqdef:Tgen}
T^{\lambda, u}:= \frac{1}{c }\Nad (\partial_x^{\rm ad} \delta_m (g + \lambda \cdot \varphi) + u ) ,
\end{equation}
where $c$ is uniquely defined by $\Vert T^{\lambda, u} \Vert_{\mathbb{L}^p (\mu) } = 1$. \footnote{see Remark \ref{rem:existence normalisation and distance} for the existence of $c$}
Define
\begin{equation}\label{eqdef:operator}
\mathcal{L}_{\varphi , \Ec }: \Theta \in \mathbb{L}^p_{\rm ad} (\mu) \mapsto \big(\E[ (\partial_x^{\rm ad} \delta_m \varphi)\Theta ], \, v \mapsto \E^\mu [ v \cdot \Theta ] \big) \in \R^k \times ( \clEconsdist )^*
\end{equation}
where $( \clEconsdist )^*$ is the topological dual of $ ( \clEconsdist , \Vert \cdot \Vert_{\mathbb{L}^{p'} (\mu ) } )$.

\noindent In this section, we prove the required differentiability result under sufficient conditions, reported in the following.
We will subsequently prove that those conditions are fulfilled for each of our situations.
\begin{AssumptionD}\label{ass:gen prop}
\begin{enumerate}[label=\textnormal{(\roman*)}]
\item \label{cond:existence argmin gen lemma} There exists $ (\hat{\lambda}, \hat{u})$ such that
$
\mathcal{L}_{\varphi , \Ec } ( T^{ \hat{\lambda},\hat{u} } ) = 0 $.
\item \label{cond:distance to constraints} There exists $C >0$ and a sequence $(\Theta_n)_n$, such that $\Theta_n \in \mathbb{L}^p_{\rm ad} (\mu)$ admits a $C^1_b$ representation and satisfies for all $n \in \N$
$$ \Vert \Theta_n - T^{ \hat{\lambda},\hat{u} } \Vert_{\mathbb{L}^p_{\rm ad} (\mu)} \rightarrow 0
\,\,
\text{and}
\,\,
\limsup_{r \rightarrow 0} \frac{1}{r} \W_p^{\rm ad}(\mu^{r \Theta_n}, \varphi^{-1}(\{0\}) \cap \Ec^{\perp} )
\leq
C\Vert \mathcal{L}_{\varphi, \Ec} ( \Theta_n ) \Vert
.
$$
\end{enumerate}
\end{AssumptionD}

\begin{Lemma}\label{lemma:gen differentiability}
Let $\varphi$ and $ g$ satisfy Assumption \ref{ass:regul and estim on g}, with $\varphi(\mu) = 0 $ and $ \mu \in \Econsgen^{\perp}$.
Assume also that Condition \hyperref[ass:gen prop]{$\textbf{\rm D}$} holds.
Then, both maps $\upGgenconsddist$ and $\lowGgenconsddist$ are differentiable at $0$, and
\begin{equation*}
\left. \upGgenconsddist\right.'  (0) = \left.\lowGgenconsddist \right.' (0)
 =
 \inf_{ \substack{ (\lambda, u) \in \R^k \times \clEconsdist } }
  \left. U^{\mathcal{C}} (\lambda, u) \right.^{1/p'}
=
 \left. U^{\mathcal{C}} (\hat{\lambda} , \hat{u})\right.^{1/p'} ,
\end{equation*}
where
 $U^{\mathcal{C}} (\lambda, u) :=
\Vert \partial_x^{\rm ad} \delta_m (g + \lambda \cdot \varphi) + u \Vert_{\mathbb{L}^{p'} (\mu)}^{p'}$.
\end{Lemma}

\noindent \textbf{Proof}.
It follows from the same line of argument as in the proof of Proposition $3.3$ of \citeauthor{Touzisauldubois2024ordermartingalemodelrisk} \cite{Touzisauldubois2024ordermartingalemodelrisk}, that
\begin{equation}\label{ineq:limsup gen case}
\begin{split}
\varliminf_{ r \rightarrow 0}
\frac{ \upGgenconsddist (r) - \upGgenconsddist (0)}{r}
\leq
\varlimsup_{ r \rightarrow 0}
\frac{ \upGgenconsddist (r) - \upGgenconsddist(0)}{r}
&\leq
\inf_{ \substack{ (\lambda, f) \in \R^k \times \Econsgen}}\left. U^{\mathcal{C}} ( \lambda,\partial_x^{\rm ad} f) \right.^{1/p'}
\\
&=
  \inf_{ \substack{ (\lambda, u ) \in \R^k \times \clEconsdist}} \left. U^{\mathcal{C}} (\lambda,u) \right.^{1/p'}
:= \ell ,
\end{split}
\end{equation}
by definition of $\clEconsdist$.
Define the following linear approximation,
$$
\hat{G}^{\mathcal{C}} (r):= \sup_{ \mu' \in \Bdisgen (\mu, r) } \hat{g}(\mu') \,\, \text{ where } \,\, \hat{g}(\mu') := g(\mu) + \int_\X \delta_m g (\mu, x ) \mu'( \mathrm{d}x)
.$$
Following the same steps as in the proof of Lemma $6.1$ of \citeauthor{Touzisauldubois2024ordermartingalemodelrisk} \cite{Touzisauldubois2024ordermartingalemodelrisk}, we obtain,
\begin{equation}
\varliminf_{ r \rightarrow 0}  \frac{ \lowGgenconsddist (r) - \lowGgenconsddist(0)}{r} = \varliminf_{ r \rightarrow 0}   \frac{ \hat{G}^{\mathcal{C}} (r) -  \hat{G}^{\mathcal{C}}(0)}{r}
.
\end{equation}
Now, by Assumption \hyperref[ass:gen prop]{$\textbf{\rm D}$} \ref{cond:existence argmin gen lemma}, there exists $ T^{ \hat{\lambda},\hat{u} } \in \mathbb{L}^p_{\rm ad} (\mu) $ such that $\mathcal{L}_{\varphi , \Ec } ( T^{ \hat{\lambda},\hat{u} } ) = 0$.
Let $(\Theta_n)_n$ be defined by Condition \ref{cond:distance to constraints} of Condition \hyperref[ass:gen prop]{$\textbf{\rm D}$}.
Let $\gamma > 0 $ and $ n \in \N$.
By Condition \ref{cond:distance to constraints}, there exists a family of measures $(\nu_{r, n})_r$ such that for all $r > 0 $, $ \nu_{r, n} \in  \varphi^{-1}(\{0\}) \cap \Ec^{\perp} $ and
\begin{equation}\label{eq:estim approx}
\W_p^{\rm ad}(\mu^{r \Theta_n}, \nu_{r, n}) \leq r ( C \Vert \mathcal{L}_{\varphi , \Ec } ( \Theta_n ) \Vert + \gamma ).
\end{equation}

Furthermore, since $\Theta_n$ is $C^1_b$, we have $ \W_p^{\rm ad}( \mu^{r \Theta_n}, \mu) \leq r \Vert \Theta_n \Vert_{\mathbb{L}^p_{\rm ad} (\mu ) }$ (consider the coupling $ \pi_r := \mu \circ ( X, X + r \Theta_n ) $, which is bi-causal since, for $r$ small enough, $\Theta_n$ is $C^1$ and $x_1 \mapsto x_1 + r (\Theta_n)_1 (x_1 ) $ is a homeomorphism), hence, by the triangle inequality,
$$
\W_p^{\rm ad}( \mu, \nu_{r, n})
\leq
r ( \Vert \Theta_n \Vert + C \Vert \mathcal{L}_{\varphi , \Ec } ( \Theta_n ) \Vert + \gamma  )
.$$
Now, $ \mathcal{L}_{\varphi , \Ec } $ is continuous and $\mathcal{L}_{\varphi , \Ec } ( T^{ \hat{\lambda},\hat{u} } ) = 0$, so, for $ n $ sufficiently large (independently of $r$), $C\Vert \mathcal{L}_{\varphi , \Ec } ( \Theta_n ) \Vert \leq \gamma$, and similarly, since $\Vert T^{ \hat{\lambda},\hat{u} }  \Vert = 1  $, $ \Vert \Theta_n \Vert \leq 1 + \gamma$.
Hence $\W_p^{\rm ad}( \mu, \nu_{r, n}) \leq r ( 1 + 3\gamma) $ and
\begin{equation}\label{ineq:limfin G hat}
\begin{split}
 (1 + 3\gamma )
\varliminf_{ r \rightarrow 0}
\frac{ \hat{G}^{\mathcal{C}} (r) - \hat{G}^{\mathcal{C}}(0)} {r}
&=
\varliminf_{ r \rightarrow 0}
\frac{ \hat{G}^{\mathcal{C}} (r (1 + 3\gamma ) ) -\hat{G}^{\mathcal{C}}(0)} {r}
\\
&\geq
\varliminf_{r \rightarrow 0} \frac{1}{r}  \int_{\X} \delta_m g(\mu, x) ( \nu_{r, n} - \mu )(\mathrm{d} x )
.
\end{split}
\end{equation}
We now prove that there exists $C >0$ such that
\begin{equation}\label{ineq:lipschitz}
\big\vert \int_{\X} \delta_m g(\mu, x) (\nu_{r, n} - \mu^{r \Theta_n }) (\mathrm{d}x )\big\vert
\leq
C \W_p^{\rm ad}( \nu_{r, n}, \mu^{r \Theta_n })
.
\end{equation}
Let $ \pi \in \Pi(\nu_{r, n}, \mu^{r \Theta_n })$.
By Assumption \ref{ass:regul and estim on g} on $g$ and applying successively the triangle inequality and Hölder's inequality,
\begin{align*}
\Big\vert \int_{\X} \delta_m g(\mu, x) (\nu_{r, n} - \mu^{r \Theta_n }) (\mathrm{d}x )\Big\vert
&=
\big\vert \E^{\pi} [ \delta_m g(\mu, X) - \delta_m g(\mu, X') ] \big\vert
\\
&\leq
\E^{\pi} \big[ \vert  X -  X'\vert^p \big]^{1/p}
\int_0^1  \E^{\pi} \big[ \vert\partial_x \delta_m g(\mu, \bar{X}_\lambda) \vert^{p'} \big]^{1/p'} \mathrm{d} \lambda
,\end{align*}
where $\bar{X}_\lambda := X' + \lambda (X - X')$.
Hence, by the Assumption \ref{ass:regul and estim on g}, $\partial_x \delta_m g$ has $(p-1)$-polynomial growth, so that $\vert \partial_x \delta_m g \vert^{p'}$ has $p$-polynomial growth, hence
$$
\int_0^1  \E^{\pi} \big[ \vert\partial_x \delta_m g(\mu, \bar{X}_\lambda) \vert^{p'} \big]^{1/p'} \mathrm{d} \lambda
\leq C ( 1 + \Vert X \Vert^{p-1}_{\mathbb{L}^{p}( \mu^{r \Theta_n} ) } +
\Vert X \Vert^{p-1}_{\mathbb{L}^{p}( \nu_{r, n} ) } )
$$
and by Inequality \eqref{eq:estim approx}, $ \Vert X \Vert_{\mathbb{L}^{p}( \nu_{r, n} ) } \leq C + \Vert X \Vert_{\mathbb{L}^{p}( \mu^{r \Theta_n} )}$, and $\Vert X \Vert_{\mathbb{L}^{p}( \mu^{r \Theta_n} )}$ is bounded, proving the desired Inequality \eqref{ineq:lipschitz} as $\pi$ is arbitrary.
Furthermore, by dominated convergence, we easily get $ \lim_{r \rightarrow 0} \frac1{r} \int_\X \delta_m g ( \mu, x) ( \mu^{r \Theta_n} - \mu )(\mathrm{d}x) = \E^{\mu} [ \partial_x \delta_m g \cdot \Theta_n ]$, hence, we have
$$
\varliminf_{r \rightarrow 0} \frac{1}{r} \int_{\X} \delta_m g(\mu, x) (\nu_{r, n} - \mu)(\mathrm{d} x)
\geq
\E^{\mu} [ \partial_x \delta_m g \cdot \Theta_n ]
-
C \W_p^{\rm ad}( \nu_{r, n}, \mu^{r \Theta_n })
\geq
\E^{\mu} [ \partial_x \delta_m g \cdot \Theta_n ]
-
2 C \gamma
.$$
Letting $\gamma$ tend to $0$ and $n$ go to infinity, by Estimate \eqref{ineq:limfin G hat},
\begin{equation}\label{ineq:liminf gen case}
\varliminf_{ r \rightarrow 0}
\frac{ \lowGgenconsddist (r)  - \lowGgenconsddist(0)}{r}
 \geq
 \E^\mu \big[ \partial_x \delta_m g \cdot T^{ \hat{\lambda},\hat{u} } \big]
\end{equation}
and $ \E^\mu \big[ \partial_x \delta_m g \cdot  T^{ \hat{\lambda},\hat{u} } \big] = \Vert \partial_x^{\rm ad} \delta_m( g + \hat{\lambda} \cdot \varphi) + \hat{u} \Vert_{\mathbb{L}^{p'} (\mu)} $
which is a consequence of Condition \hyperref[ass:gen prop]{$\textbf{\rm D}$} \ref{cond:existence argmin gen lemma}.
Indeed $ \Vert \partial_x^{\rm ad} \delta_m( g + \hat{\lambda} \cdot \varphi) + \hat{u} \Vert_{\mathbb{L}^{p'} (\mu)}
 =
 \E^\mu [ (\partial_x^{\rm ad} \delta_m( g + \hat{\lambda} \cdot \varphi) + \hat{u} )\cdot T^{ \hat{\lambda},\hat{u} } ]
 $ and by Condition \hyperref[ass:gen prop]{$\textbf{\rm D}$} \ref{cond:existence argmin gen lemma} applied to $v = \hat{u}$,
 $$
\E^\mu [ (\partial_x^{\rm ad} \delta_m (\hat{\lambda} \cdot \varphi) ) \cdot T^{ \hat{\lambda},\hat{u} } ] = 0
\,\, \text{and}\,\,
\E^\mu [ \hat{u} \cdot T^{ \hat{\lambda},\hat{u} } ] = 0
.$$
Finally, noting that $ \Vert \partial_x^{\rm ad} \delta_m( g + \hat{\lambda} \cdot \varphi) + \hat{u} \Vert_{\mathbb{L}^{p'} (\mu)} \geq \ell$, we proved that $ \varliminf_{ r \rightarrow 0}
\frac{ \lowGgenconsddist (r)  - \lowGgenconsddist(0)}{r}  \geq \ell$.
Hence, since $ \lowGgenconsddist \leq \upGgenconsddist$, putting Equations \eqref{ineq:limsup gen case} and \eqref{ineq:liminf gen case} together proves the differentiability of $ \upGgenconsddist $ and $ \lowGgenconsddist $ at $0$, with both derivatives equal to $ \ell = \left. U^{\mathcal{C}} (\hat{\lambda} , \hat{u})\right.^{1/p'} $.
\ep

\subsection{The martingale Coupling Case}

For the martingale coupling case, the proof relies on the sufficient condition established in Subsection~\ref{subsection:sufficientcondition}.
However, the main difficulty lies in constructing a family of martingale couplings $\nu_r$ that are a good approximation to $\text{argmax}_{ \mu' \in \BadMm (\mu, r) } g(\mu')$. 
The construction of such a family of measures relies on an adaptation of the Implicit Function Theorem; see Subsection \ref{subsec:implicit function theorem}.

We introduce the following notation.
Given a family of random variables $ (X_n)_n $, where for $ n \in \N$, $X_n $ is defined on a probability space $(\Omega, \Fc, \P) $, we say that
\begin{equation}\label{eqdef:convergencerandomvar}
\text{$ X_n = o_\P (1)$ if $ X_n \xrightarrow[n \rightarrow \infty]{ \P} 0 $ and $ X_n =  o_{\mathbb{L}^p (\P)} (1)$ if $\Vert X_n \Vert_{\mathbb{L}^p (\P)} \xrightarrow[n \rightarrow \infty]{} 0 $.}
\end{equation}
For a sequence $r_n \rightarrow 0$, we write $X_n = o_\P (r_n)$ (resp. $X_n = o_{\mathbb{L}^p (\P)} (r_n)$) if $X_n / r_n = o_\P (1)$ (resp. $X_n / r_n = o_{\mathbb{L}^p (\P)} (1)$), and we use the same notations with $\P$ replaced by $\mu$ or $\nu$ when the random variables are defined on the canonical space.

\noindent In the following, for $(E,\vert \cdot \vert_E)$ a normed vector space, we let
\begin{equation}\label{eqdef:continuousoperator}
\Lc_c(E) := \{ \Psi:(E,\vert \cdot \vert_E)\to(E,\vert \cdot \vert_E) \,\,\ \text{linear and continuous} \} ,
\end{equation}
 and we denote by $\vertiii{\Psi}$ the corresponding operator norm
\begin{equation}\label{eqdef:def operator norm}
    \vertiii{\Psi} := \sup_{\vert x \vert_E \leq 1} \vert \Psi(x) \vert_E .
\end{equation}

\noindent Define the operator $ \E^{\mu}_1 \circ \E^{\mu}_2 [ f ] = \E^{\mu} [ \E^{\mu} [ f(X_1) \vert X_2 ] \vert X_1 ] $ for $f \in \mathbb{L}^{\alpha}(\sigma(X_1), \mu)$. 
We note that the operator $\E^{\mu}_1 \circ \E^{\mu}_2 $ maps $ \mathbb{L}^{\alpha}_0 (\mu_1)$ to $\mathbb{L}^{\alpha}_0 (\mu_1) $, see definition \eqref{eqdef:L0} of $\mathbb{L}^{\alpha}_0(\mu_1)$.

\begin{Lemma}\label{lemma:closedness ans coercivity}
Let Assumption \hyperref[ass:mart coupling]{$\textbf{\rm C}_{\W_p^{\rm ad}}$} hold and $ 1 \leq \alpha \leq \infty$. 

\begin{enumerate}[label=\textnormal{(\roman*)}]
\item \label{contraction_operator} Then $ \vertiii{ \E^{\mu}_1 \circ \E^{\mu}_2}_\alpha < 1 $.
\item \label{direct_sum_sumspace} The vector space $ \Ec_\alpha:= \{ f_1 \oplus f_2 , f_1 \in \mathbb{L}^{\alpha}_0 (\mu_1) , \,\, f_2 \in \mathbb{L}^{\alpha} (\mu_2) \} $ can be written as the direct sum $ \Ec_\alpha = \mathbb{L}^{\alpha}_0 (\mu_1) \oplus \mathbb{L}^{\alpha} (\mu_2) $ and is closed in $\mathbb{L}^{\alpha} (\mu) $.  
\item \label{coercivity direct sum} There exists $ C > 0 $ such that for all $ f_1 \oplus f_2 \in \Ec_\alpha $, we have 
$$
\Vert f_1 \oplus f_2 \Vert_{\mathbb{L}^{\alpha} (\mu)} \geq C(\Vert f_1 \Vert_{\mathbb{L}^{\alpha} (\mu_1)} + \Vert f_2 \Vert_{\mathbb{L}^{\alpha} (\mu_2)})
.$$
\end{enumerate}
\end{Lemma}

\noindent The proof of this lemma is deferred to the end of this section. 

\begin{Lemma}\label{lemma:existence of martingale coupling}
Let $\Theta_2: \R^2 \rightarrow \R $ be a $C^2$, compactly supported function, let $1 < \alpha < \infty$, let $\mu$ satisfy Assumption \hyperref[ass:mart coupling]{$\textbf{\rm C}_{\W_p^{\rm ad}}$}, and define
\begin{equation}\label{eqdef:hatX,Tc,F^a}
\Tc_{r, a} := F^{-1}_2 \circ F_2^{r, a} \,\, , \,\, 
F_2^{r,a} (x):= \mu (\hat{X}_2^{r, a}\leq x ) 
\,\, \text{where} \,\,
\hat{X}_2^{r, a}:= X_2 + r \Theta_2 + a(X_1) 
.\end{equation} 
Then, there exists $\eta > 0 $ and $a: (-\eta, \eta) \rightarrow \mathbb{L}^p_0 (\mu_1)$ such that:
$$
\E^{\mu}_1 [ \Tc_{r, a_r} (\hat{X}_2^{r, a_r}) ] = X_1
\,\,
\text{for all $ \vert r \vert < \eta $.} 
$$
Moreover, for some constant $C >0$,
\begin{equation} \label{ineq:estim ar}
\Vert a_r \Vert_{\mathbb{L}^{\alpha} (\mu_1)} \leq r C (\Vert \E^{\mu}_1 [ \Theta_2 - \E^{\mu}_2 [ \Theta_2 ] ] \Vert_{\mathbb{L}^{\alpha} (\mu_1)} + o (1) ) 
.\end{equation}
\end{Lemma}

\noindent The proof of this lemma is deferred to later.

\begin{Lemma}\label{lemma:proj mart coupling set}
Let $\mu$ satisfy Assumption \hyperref[ass:mart coupling]{$\textbf{\rm C}_{\mathbb{W}^{\rm ad}_p}$}.
Let $ \Theta_2: \R^2 \rightarrow \R$ be a compactly supported $C^2$ function, and set $\Theta := (0, \Theta_2)$.
Then there exists a family of measures $ (\nu_r)_{ r > 0 }$ such that
\begin{equation}
 \nu_r \in \Pi^{\rm M} (\mu_1, \mu_2)
\,\, \text{and} \,\,
\W_p^{\rm ad} (\mu^{r\Theta}, \nu_r) \leq C r \big(\Vert \opedistmartcoup (\Theta) \Vert + o (1) \big) \,\, \text{for all } \,\, r > 0
.
\end{equation}
\end{Lemma}

\noindent We also postpone the proof of this lemma and move on to the proof of Propositions \ref{prop:sensi 1 dim mart coupling} and \ref{prop:proj mart coupling set}.

\medskip

\noindent \textbf{Proof of Proposition \ref{prop:proj mart coupling set}.} 
By Lemma \ref{lemma:proj mart coupling set}, there exists a family $\nu_r \in \Pi^{\rm M} (\mu_1, \mu_2)$, for which we have 
$$
\frac{1}r \W_p^{\rm ad}( \mu^{r \Theta}, \Pi^{\rm M} (\mu_1, \mu_2) ) 
\leq 
\frac{1}{r}
\W_p^{\rm ad} ( \mu^{r \Theta}, \nu_r ) 
\leq 
 C (\Vert \opedistmartcoup (\Theta) \Vert + o (1) )
.$$ 
Letting $r$ tend to $0$, we get the desired result.
\ep 

\medskip 

\noindent \textbf{Proof of Proposition \ref{prop:sensi 1 dim mart coupling}.}
We rely on Lemma \ref{lemma:gen differentiability} of Subsection~\ref{subsection:sufficientcondition}, and proceed in two steps.

\noindent{\bf{Step $1$.}} 
We first prove that Condition \hyperref[ass:gen prop]{$\textbf{\rm D}$} \ref{cond:existence argmin gen lemma} is satisfied. 
Using the notation \eqref{eqdef:oplusnotation}, define $ \Econsgen_{\rm M, m}:= \{ f_1 \oplus f_2 - \mu (f_1 \oplus f_2) + h^{\otimes} , f_1, f_2, h \in C^1_b(\R, \R) \} .$ 
In the sense of \eqref{eqdef:orthogonalityduality}, we have $\Pi^{\rm M} (\mu_1, \mu_2) = \Econsgen_{\rm M, m}^{\perp}$.
Let $ u = f_1 \oplus f_2 -\mu (f_1 \oplus f_2) + h^{\otimes} $ in $\Econsgen_{\rm M, m}$; since $\mu$ is a martingale measure,
$
\partial_x^{\rm ad} u (X) = (f_1', f_2')^{\intercal} + h J$, with $J$ as in Proposition \ref{prop:sensi 1 dim mart coupling}.
As $C^1_b$ is dense in $\mathbb{L}^{p'} (\mu_i)$ for $i=1,2$, we can extend this result to its closure in $\mathbb{L}^{p'}$.
Hence, $\clEconsdist_{\rm M, m}$, defined by \eqref{eqdef:closure image of gradient}, is
$$
\clEconsdist_{\rm M, m} = \mathbb{L}^{p'} (\mu_1) \times ( \mathbb{L}^{p'} (\mu_1) + \mathbb{L}^{p'} (\mu_2) )
= \mathbb{L}^{p'} (\mu_1) \times ( \mathbb{L}^{p'}_0 (\mu_1) \oplus \mathbb{L}^{p'} (\mu_2) ),
$$
where the second equality follows directly from Lemma \ref{lemma:closedness ans coercivity}. 
Let $ U_{\rm ad}^{\rm M, m} $ be defined as in Proposition \ref{prop:sensi 1 dim mart coupling}. 
By Lemma \ref{lemma:closedness ans coercivity}, $ U_{\rm ad}^{\rm M, m} $ is coercive (in the sense of Definition \ref{def:coercivity}). 
Furthermore, since $p' > 1$ and the sum $\mathbb{L}^{p'}_0 (\mu_1) \oplus \mathbb{L}^{p'} (\mu_2)$ is direct, $ U_{\rm ad}^{\rm M, m} $ is strictly convex.
Hence $ U_{\rm ad}^{\rm M, m} $ is continuous, coercive, strictly convex, and $\clEconsdist_{\rm M, m} $ is a reflexive Banach space, being a closed subspace of a reflexive Banach space. 
Therefore, there exists a minimizer which satisfies the first-order condition given by \eqref{eqdef:FOC Mart couplings}. 
A direct consequence of the first-order condition \eqref{eqdef:FOC Mart couplings} is that the minimizer $(f^{\rm M, m}, h_{\rm M, m})$ satisfies Condition \hyperref[ass:gen prop]{$\textbf{\rm D}$} \ref{cond:existence argmin gen lemma}.

\smallskip

\noindent{\bf{Step $2$.}} 
Condition \hyperref[ass:gen prop]{$\textbf{\rm D}$} \ref{cond:distance to constraints} is satisfied as a direct consequence of Proposition \ref{prop:proj mart coupling set}.
\smallskip
\noindent Applying Lemma \ref{lemma:gen differentiability} yields the desired result.
\ep 

\medskip

\noindent We now move on to the proof of technical lemmas. Most of them rely on the expansion Lemma \ref{lemma:expansion in Lp cdf}, whose statement and proof are deferred to Subsection \ref{Expansions in Probability and}.

\medskip

\noindent\textbf{Proof of Lemma \ref{lemma:proj mart coupling set} }
Let $\Theta_2: \R^2 \rightarrow \R $ be a compactly supported $C^2$ function, and set $\Theta = ( 0, \Theta_2)$. 
Let $\eta >0$ and $a: (-\eta, \eta) \rightarrow \mathbb{L}^p_0 (\mu_1)$ be the map constructed in Lemma \ref{lemma:existence of martingale coupling}. 
Using Notations \eqref{eqdef:hatX,Tc,F^a}, we let $X'_1 = X_1$ and $X'_2 = \Tc_{r, a_r } ( \hat{X}^{r, a_r}_2 )$.
By construction, we have $\nu_r:= \mathcal{L} (X') \in \Pi (\mu_1, \mu_2)$.
Let $\Theta_2 + \frac{1}{r} a_r(X_1)$, is $p$-uniformly integrable since $ (  \frac{a_r}{r} )_r $ is bounded in $ \mathbb{L}^\alpha (\mu_1)$, and $\alpha$ can be chosen strictly greater than $p$ in Lemma \ref{lemma:existence of martingale coupling}. 
Therefore, by Lemma \ref{lemma:expansion in Lp cdf} applied to $\mu$ and $\Theta + \frac{1}{r} a_r$, we obtain
\begin{equation*}
\begin{split}
F_2^{ r, a_r} (\hat{X}_2^{r, a_r}) 
=
F_2 (X_2) &+ r\int_\R (\Theta_2 (X) - \Theta_2(x_1, X_2) ) \kappa(x_1, X_2) \mu_1 (\mathrm{d}x_1) 
\\
&+ \int_\R (  a_r (X_1) - a_r (x_1)  ) \kappa(x_1, X_2) \mu_1 (\mathrm{d}x_1)  + o_{\mathbb{L}^p (\mu)} (r)
.\end{split}
\end{equation*}
Recall from Assumption \hyperref[ass:mart coupling]{ $ \textbf{\rm C}_{\W_p^{\rm ad} } $} \ref{cond:density second marginal} that $F_2^{-1} \in C^1 ((0, 1) , I ) $ and, $ q_2 (X_2) \geq c > 0$, $\mu_2$ almost surely, which ensures that $F_2^{-1}$ is uniformly Lipschitz from $(0, 1)$ into $I$. 
Therefore, by Lemma \ref{lemma:expans in Lp C1},  
\begin{equation}\label{eq:expansion second component}
X'_2
=
X_2 + r(\Theta_2 - \E^{\mu}_2 [\Theta_2]) + (a_r(X_1) - \E^{\mu}_2 [ a_r (X_1) ]) + o_{\mathbb{L}^p (\mu)} (r)
.\end{equation}
Now consider the coupling $ \pi := \Lc (X + r\Theta(X), X') $. It is clearly bi-causal since $X_1 = X'_1$. 
Furthermore, 
$$
\frac{1}{r}\W_p^{\rm ad}( \mu^{r \Theta}, \Pi^{\rm M}(\mu_1, \mu_2) ) \leq \frac{1}{r}\E^{\pi} [ \vert X - X' \vert^p]^{1/p}
.$$
Moreover, using the expansion \eqref{eq:expansion second component} and the estimate \eqref{ineq:estim ar} from Lemma \ref{lemma:existence of martingale coupling}, we obtain the desired result.
\ep

\medskip

\noindent We now move on to the proof of Lemma \ref{lemma:existence of martingale coupling}.
Let $ \Theta_2 : \X \rightarrow S$ be a $C^2$, compactly supported function. Let $ 1 < \alpha < \infty$ and define
 \begin{equation}\label{eqdef:operateur}
 \begin{split}
  \Uc &: \R \times \mathbb{L}^\alpha_0 (\mu_1) \longrightarrow \mathbb{L}^\infty_0 (\mu_1)\\
  & (r, a) \mapsto \E^{\mu}_1 [ \Tc_{r, a } ( \hat{X}^{r, a}_2 ) ] - X_1 
\end{split}
.\end{equation}

\begin{Lemma}\label{lemma:operateur is gateau diff} 
Assume that \hyperref[ass:mart coupling]{$\textbf{\rm C}_{\W_p^{\rm ad}}$} holds. 
Then $\Uc$ is well defined and Gâteaux differentiable in a convex neighborhood $I \times \mathcal{V} $ of $(0, 0)$, with Gâteaux derivative given by 
\begin{equation}\label{eqdef:gateau derivative}
 \begin{split}
\partial_r \Uc (r, a) 
 &:=
 \int_\R
 \E^{\mu}_1 \Big[ 
 (\hat{\Theta}(\hat{X}) - \hat{\Theta} (x_1, \hat{X}_2^{r, a}) )
\frac{ \kappa_{r, a}(x_1, \hat{X}_2^{r, a})}{q_2 \circ \Tc_{r, a} (\hat{X}_2^{r, a}) 
} 
\Big]
 \mu_1(\mathrm{d} x_1) 
\\
\mathrm{d}_a \Uc (r, a) (b) 
 &:=
 \int_\R
 \E^{\mu}_1 \Big[
  (b(X_1) - b(x_1) )
 \frac{\kappa_{r, a}(x_1, \hat{X}_2^{r, a})
}{ q_2 \circ \Tc_{r, a} (
 \hat{X}_2^{r, a}  )
} 
\Big]
 \mu_1(\mathrm{d} x_1)
,\end{split}
 \end{equation}
 where we used notations \eqref{eqdef:hatX,Tc,F^a}, $ \hat{\Theta} (x_1, x_2):= \Theta_2 (x_1, \varphi (r, x_1, x_2 - a(x_1) ) )$, $ \hat{X}_1 = X_1$,  $ \varphi: \R \times \R^2 \rightarrow \R $ is the inverse of $ x_2 \mapsto x_2 + r \Theta_2 (x) $ and 
\begin{equation}\label{eqdef:kappatransition}
\kappa_{r, a} (x) := 
\kappa (x_1, \varphi (r, x_1, x_2 - a(x_1) ) ) \partial_{x_2} \varphi (r, x_1, x_2 - a(x_1) ) ,
\end{equation}
is the disintegration of $\mu\circ (X_1, \hat{X}_2^{r, a})^{-1}$ with respect to its first marginal. 
 \noindent Furthermore,  for $(r, a) \in  I \times \mathcal{V} $ and $\beta \in \{ \alpha, \infty\}$, $\mathrm{d}_a \Uc (r, a) \in \mathcal{L}_c (\mathbb{L}^{\beta}_0(\mu_1) )$.
 \end{Lemma}

\noindent\textbf{Proof of Lemma \ref{lemma:operateur is gateau diff}.}
We first check that $\Uc$ is well defined. 
For any $r > 0$ and $a \in \mathbb{L}^{\alpha}_0 (\mu_1)$, we have $ \Tc_{r, a} (\hat{X}_2^{r, a}) \sim \mu_2 $. 
Since $ F_2^{-1} $ is bounded by Assumption \hyperref[ass:mart coupling]{ $ \textbf{\rm C}_{\W_p^{\rm ad} } $} \ref{cond:density second marginal} and $X_1$ is essentially bounded, we have $ \Uc (r, a) \in \mathbb{L}^\infty (\mu_1)$ and $ \E^{\mu}[ \Uc(r, a) ] = \int_\R x \mu_2(\mathrm{d}x) - \int_\R x \mu_1(\mathrm{d}x) = 0 $.

\smallskip

\noindent We now establish the Gâteaux differentiability of $\Uc$ with respect to $(r, a)$. 
Let $ a \in \mathbb{L}^\alpha_0 (\mu_1) $, $r$, $\varepsilon:= (\varepsilon_1, \varepsilon_2) >0$, small enough.
Define $\hat{X}_1:= X_1$, and let $\hat{X}_2^{r, a}$, $ \hat{\Theta} (x_1, x_2)$, and $\varphi(r, x_1, y_2)$ be as in \eqref{eqdef:gateau derivative}. 
The measure $ \nu_{r, a}:= \mathcal{L} (\hat{X}) $ admits the disintegration 
$
\nu_{r, a} = \kappa_{r, a} (x) \mathrm{d}x_2 \mu_1 (\mathrm{d}x_1)
$ where $\kappa_{r, a}$ defined by \eqref{eqdef:kappatransition}.
Applying Lemma \ref{lemma:expansion in Lp cdf} to the measure $\nu_{r, a}$ with $\tilde{\Theta}_\varepsilon := \frac{\varepsilon_1}{\vert \varepsilon \vert } \hat{\Theta} (\hat{X}) + \frac{\varepsilon_2}{\vert \varepsilon \vert} b(\hat{X}_1)$ (which is $p$-uniformly integrable), yields
\begin{align*}
F_2^{ r + \vert \varepsilon \vert \tilde{\Theta}_\varepsilon} \big(\hat{X}_2^{r, a} + \vert \varepsilon \vert \tilde{\Theta}_\varepsilon ( \hat{X} ) \big) 
&=
F_2^{r,a} (\hat{X}_2^{r, a})
+
\vert \varepsilon \vert  \mathcal{L}_{ r, a} (\tilde{\Theta}_\varepsilon) + o_{\mathbb{L}^p (\mu)} (\varepsilon ) 
\\
& 
=
F_2^{r,a} (\hat{X}_2^{r, a})
+
\varepsilon_1  \mathcal{L}_{ r, a} (\hat{\Theta}) + \varepsilon_2 \mathcal{L}_{ r, a} (b)
+ o_{\mathbb{L}^p (\mu)} (\varepsilon ),
 \end{align*}
where $
\mathcal{L}_{ r, a} (\hat{\Theta}):= 
 \int_\R
 (\hat{\Theta}(\hat{X}) - \hat{\Theta} (x_1, \hat{X}_2^{r, a}) ) \kappa_{r, a}(x_1, \hat{X}_2^{r, a}) \mu_1(\mathrm{d} x_1)$. 
 We now prove that  $ \mathrm{d}_a \Uc (r, a): \mathbb{L}^\alpha_0 \rightarrow \mathbb{L}^\alpha_0$ and $ \mathrm{d}_a \Uc (r, a): \mathbb{L}^\infty_0 \rightarrow \mathbb{L}^\infty_0$ are bounded linear continuous operators. 
 By Assumption \hyperref[ass:mart coupling]{ $ \textbf{\rm C}_{\W_p^{\rm ad} } $} \ref{cond:disintegration}-\ref{cond:density second marginal}, and since $\hat{\Theta}$ is bounded, we have the following inequalities:
\begin{equation}\label{eq:estim gateau mart}
\begin{split}
&\big\vert 
  \big(b(X_1) - b(x_1) \big)
 \frac{\kappa_{r, a}(x_1, \hat{X}_2^{r, a}) 
}{q_2 \circ \Tc_{r, a} \big(\hat{X}_2^{r, a} \big) 
} 
\big\vert 
\leq 
C \vert b(X_1) - b(x_1) \vert 
.\end{split}
\end{equation}
Estimates \eqref{eq:estim gateau mart} yield the desired continuity properties.
\ep 

\medskip

\begin{Lemma}\label{lemma:continuity gateau diff}
Let \hyperref[ass:mart coupling]{$\textbf{\rm C}_{\W_p^{\rm ad}}$} hold and let $I \times \mathcal{V}$ and $\partial_r \Uc $,  $\mathrm{d}_a \Uc $ be defined by Lemma \ref{lemma:operateur is gateau diff}.
Then $\partial_r \Uc $ is continuous with respect to $r$.
Let $a, \tilde{a}$ in $\mathcal{V}$ and $r \in I$, then the mapping 
$$
\lambda \in [0, 1] \mapsto \mathrm{d}_a \Uc (r, \bar{a}_\lambda) (a - \tilde{a}) \in \mathbb{L}_0^{\alpha} (\mu_1) \,\, \text{where} \,\, \bar{a}_\lambda := (1- \lambda) a + \lambda \tilde a,
$$
is continuous. 
Furthermore, the mapping $ (r, a) \in \R \times \mathbb{L}^{\infty}_0 (\mu_1) \mapsto \mathrm{d}_a \Uc (r, a) \in \mathcal{L}_c (\mathbb{L}^{\beta}_0 (\mu_1) ) $ is continuous at $(0, 0)$, for $ \beta \in \{ \alpha, +\infty\}$.
\end{Lemma}

\noindent \textbf{Proof of Lemma \ref{lemma:continuity gateau diff}.} 
We begin with the first point. 
Let $a $ and $ \tilde{a}$ in $ \mathbb{L}^\alpha_0 (\mu_1) $, we will prove the continuity of $\lambda \mapsto \Lambda (\lambda) \in \mathbb{L}^\alpha (\mu_1) $, where $\Lambda (\lambda):= \mathrm{d}_a \Uc (r, \bar{a}_\lambda ) (a - \tilde{a}) $, with $\bar{a}_\lambda$ as in the statement.
Without loss of generality, we only check the continuity at $\lambda = 0$. 
Since $q_2 \circ F_2^{-1} $ and $\kappa$ are bounded, we have 
$$
\vert 
\Lambda(\lambda) - \Lambda(0) 
\vert 
\leq I_1 (\lambda) + I_2 (\lambda), $$
where 
$$
I_1 (\lambda)
:=
\int_{\R} \E_1^{\mu} \Big[ k(x_1, X_1) \big\vert  \kappa_{r, \bar{a}_\lambda} (x_1, \hat{X}_2^{r, \bar{a}_\lambda} ) - 
\kappa_{r, a} (x_1, \hat{X}_2^{r, a})  
\big\vert \Big]\mu_1 (\mathrm{d}x_1),
$$
with $k(x_1, X_1) = \vert a - \tilde{a} \vert (x_1) + \vert a - \tilde{a} \vert (X_1)$ and  
$$
I_2 (\lambda)
:=
\int_{\R} \E_1^{\mu} \Big[ \big\vert q_2 \circ \Tc_{r, \bar{a}_\lambda} ( \hat{X}_2^{r, \bar{a}_\lambda}) -  q_2 \circ \Tc_{r, a} ( \hat{X}_2^{r, a}) 
\big\vert \Big]\mu_1 (\mathrm{d}x_1)
.$$
By Lemma \ref{lemma:convergence of integrals}, $\Vert I_1 (\lambda) \Vert_{\mathbb{L}^p} \rightarrow 0$. 
Furthermore, since for all $0\leq \lambda \leq 1$, $\mu \circ ( X_2 + r \Theta_2(X) + \bar{a}_\lambda (X_1) )^{-1}$ admits a density, we have, as $\lambda$ tends to $0$ that $F^{r, \bar{a}_\lambda}_2 \rightarrow F^{r, a}_2 $ uniformly. 
Since $q_2 \circ F_2^{-1}$ is bounded and uniformly continuous, standard arguments yield $\Vert I_2 (\lambda) \Vert_{\mathbb{L}^p (\mu)} \rightarrow 0$.

\noindent We now prove that $ (r, a) \in (-\eta, \eta) \times \mathbb{L}^{\infty}_0 (\mu_1) \mapsto \mathrm{d}_a \Uc (r, a) \in \mathcal{L}_c (\mathbb{L}^{\beta}_0) $ is continuous at $(0, 0)$ for $\beta \in \{\alpha,  \infty \}$.
To do so, it is sufficient to prove that $ (r, a) \in (-\eta, \eta) \times \mathbb{L}^{\infty}_0 (\mu_1) \mapsto \mathrm{d}_a \Uc (r, a) \in \mathcal{L}_c (\mathbb{L}^{\beta}_0) $ is continuous at $(0, 0)$ for $1 < \beta < \infty $.
Let $ \Vert a_n \Vert_{\mathbb{L}^{\infty} (\mu) } \rightarrow 0 $, $ r_n \rightarrow 0 $, and consider an arbitrary $b \in \mathbb{L}^\beta_0 (\mu_1)$.
Then, letting $\beta' $ be the conjugate exponent of $\beta$, by Hölder's inequality, we get 
\begin{align*}
&\Vert 
\mathrm{d}_a \Uc (r_n, a_n) (b) - \mathrm{d}_a \Uc (0, 0) (b) 
\Vert_{\mathbb{L}^\beta(\mu_1)}
\\
=
&\Big\Vert 
 \int_\R
 \E^{\mu}_1 \big[
  (b(X_1) - b(x_1) )\Psi_n (x_1, X)
\big]
 \mu_1(\mathrm{d} x_1)
 \Big\Vert_{\mathbb{L}^\beta(\mu_1)}
\\
\leq 
&\Big\Vert 
\Big(
 \int_\R
 \E^{\mu}_1 \big[
  \vert b(X_1) - b(x_1) \vert^{\beta}
  \big]
   \mu_1(\mathrm{d} x_1)
   \Big)^{1/\beta}
   \Big(
 \int_\R
 \E^{\mu}_1 \big[
\vert \Psi_n (x_1, X)\vert^{\beta'}
\big]
 \mu_1(\mathrm{d} x_1)
 \Big)^{1/\beta'}
 \Big\Vert_{\mathbb{L}^\beta(\mu_1)}
\\
\leq &2 c \Vert b \Vert_{\mathbb{L}^\beta (\mu_1) }  
\Big\Vert \int_\R \E_1^\mu [ \vert \Psi_n (x_1, X) \vert^{\beta'} ] \mu_1( \mathrm{d}x_1 )
\Big\Vert_{ \mathbb{L}^{\infty}(\mu_1)}
,\end{align*}
where $\Psi_n (x_1, X):= \frac{\kappa_{r_n, a_n}(x_1, \hat{X}_2^{r_n, a_n}) 
}{ q_2 \circ \Tc_{r_n, a_n}(\hat{X}_2^{r_n, a_n}) 
}
-
\frac{\kappa(x_1, X_2) }{q_2 (X_2) }
$. 
Now, since $ q_2 \circ F_2^{-1} \geq c$ by Assumption \hyperref[ass:mart coupling]{$\textbf{\rm C}_{\W_p^{\rm ad}}$} \ref{cond:density second marginal}, $\Psi_n$ is bounded, letting $C = \Vert \Psi_n \Vert_{\mathbb{L}^\infty (\mu_1)}$, we have
\begin{align*}
\Vert 
\mathrm{d}_a \Uc (r_n, a_n) (b) - \mathrm{d}_a \Uc (0, 0) (b) 
\Vert_{\mathbb{L}^\beta(\mu_1)}
&\leq 2 c C^{\beta' - 1} \Vert b \Vert_{\mathbb{L}^\beta (\mu_1) }  
\Big\Vert \int_\R \E_1^\mu [ \vert \Psi_n (x_1, X) \vert ] \mu_1( \mathrm{d}x_1 )
\Big\Vert_{ \mathbb{L}^{\infty}(\mu_1)},
\end{align*}
which is possible since $ 1 < \beta \leq \infty $ hence $\beta' < \infty $.
Now, $ \mu_1 \otimes \mu$ almost surely 
$$
\vert \Psi_n (x_1, X) \vert \leq 
C ( I^1_n (x_1, X) + I^2_n (x_1, X) )
,$$
where $I^1_n (x_1, X) :=
\vert \kappa_{r_n, a_n}(x_1, \hat{X}_2^{r_n, a_n}) - \kappa (x_1, X_2) \vert $ and $I_n^2(x_1, X):= 
\big\vert q_2 \circ \Tc_{r_n, a_n} (\hat{X}_2^{r_n, a_n}) - q_2 \big(X_2 \big) \big\vert$. 
Note that 
$$
\int_\R \E_1^\mu [ \vert I_n^1(x_1, X) \vert ] \mu_1( \mathrm{d}x_1 )
= 
\int_\R \int_\R
\kappa(X_1, x_2) 
\vert \kappa_{r_n, a_n}(x_1, \hat{X}_2^{r_n, a_n}) - \kappa (x_1, X_2) \vert \mathrm{d}x_2 \mu_1(\mathrm{d}x_1)
.$$
Now, since $ a_n \xrightarrow[n \rightarrow \infty]{\mathbb{L}^\infty (\mu_1) } 0 $, we may apply Lemma \ref{lemma:convergence of integrals} and obtain that 
$$
\Big\Vert \int_\R \E_1^{\mu} \big[ I_n^1 (x_1, X) \big] \mu_1( \mathrm{d}x_1 ) \Big\Vert_{\mathbb{L}^\infty (\mu_1) } \rightarrow 0
.$$
Moreover,
$$
\Big\Vert \int_\R \E_1^{\mu} \big[ I_n^2 (x_1, X) \big] \mu_1( \mathrm{d}x_1 ) \Big\Vert _{\mathbb{L}^\infty (\mu_1) }
\leq 
R_n^1 + R_n^2 
,$$
where, setting $E_n := \big\vert F_2^{r_n,a_n} (\hat{X}_2^{r_n, a_n}) - F_2 (X_2) \big\vert$,
\begin{equation*}
\begin{split}
R_n^1 &:= \Big\Vert \int_\R \E^{\mu}_1 \big[ I_n^2 (x_1, X) \mathds{1}_{ \{ E_n \geq \delta \}} \big] \mu_1 (\mathrm{d}x_1) \Big\Vert_{\mathbb{L}^\infty (\mu_1)} ,
 \\
R_n^2 &:= \Big\Vert \int_\R \E^{\mu}_1 \big[ I_n^2 (x_1, X) \mathds{1}_{ \{ E_n \leq \delta \}} \big] \mu_1 (\mathrm{d}x_1) \Big\Vert_{\mathbb{L}^\infty (\mu_1)} ,
\end{split}
\end{equation*}
and where, given $\varepsilon > 0$, $\delta > 0$ is chosen, by uniform continuity of $q_2 \circ F_2^{-1}$, such that $\vert q_2 \circ F_2^{-1} (u) - q_2 \circ F_2^{-1} (v) \vert \leq \varepsilon$ whenever $\vert u - v \vert \leq \delta$.

\noindent By uniform continuity of $q_2 \circ F_2^{-1}$, and since $ F_2^{-1} (F_2 (X_2)) = X_2 $, $\mu_2$ almost surely, we have $ R_n^2 \leq \varepsilon $. 
Finally, since $q_2 $ is bounded, 
$$
R_n^1 \leq C
\Big\Vert \int_\R \E^{\mu}_1 \big[ \mathds{1}_{ \{ E_n \geq \delta \}} \big] \mu_1 (\mathrm{d}x_1) \Big\Vert_{\mathbb{L}^\infty (\mu_1)},
$$
and, by the conditional Markov inequality, we get 
$$
R_n^1 \leq \frac{C}{\delta} \Big\Vert \int_\R \E^{\mu}_1 [ E_n ] \mu_1 (\mathrm{d}x_1) \Big\Vert_{\mathbb{L}^\infty (\mu_1)} 
.$$
Now, since $ F_2 $ is continuous, as $\mu_2$ admits a density, we clearly have $ F_2^{r_n,a_n} \xrightarrow[n \rightarrow \infty]{} F_2 $ uniformly, hence 
\begin{align*}
R_n^1 &\leq \frac{C}{\delta} \Vert F_2 - F_2^{r_n, a_n} \Vert_{\mathbb{L}^{\infty} (\mu)}
\\
&\quad + \frac{C}{\delta} \Big\Vert \int_\R \E^{\mu}_1 \Big[ \big\vert F_2 (\hat{X}_2^{r_n, a_n}) - F_2 (X_2) \big\vert \Big] \mu_1 (\mathrm{d}x_1) \Big\Vert_{\mathbb{L}^\infty (\mu_1)},
\end{align*}
 and, as $ F_2 $ is Lipschitz since $q_2$ is bounded (by boundedness of $\kappa$), we obtain, up to a change in the constant $C$, 
 \begin{align*}
R_n^1 &\leq \frac{C}{\delta} \Vert F_2 - F_2^{r_n, a_n} \Vert_{\mathbb{L}^{\infty} (\mu)} + \frac{C}{\delta} \Big\Vert \int_\R \E^{\mu}_1 \Big[ 
 \big\vert \hat{X}_2^{r_n, a_n} - X_2 \big\vert \Big]
 \mu_1 (\mathrm{d}x_1) \Big\Vert_{\mathbb{L}^\infty (\mu_1)} 
 \\
 &\leq 
 \frac{C}{\delta} \Vert F_2 - F_2^{r_n, a_n} \Vert_{\mathbb{L}^{\infty} (\mu)} + \frac{C}{\delta} (r_n \Vert \Theta_2 \Vert_{\mathbb{L}^{\infty} (\mu)}+ \Vert a_n \Vert_{\mathbb{L}^\infty (\mu) }).
 \end{align*}
This proves that $ \vertiii{ \mathrm{d}_a \Uc (r_n, a_n) - \mathrm{d}_a \Uc (0, 0)}_{\mathbb{L}^\infty} \rightarrow 0 $ as $n$ goes to infinity. \ep

\medskip 

\noindent From Lemmas \ref{lemma:continuity gateau diff} and \ref{lemma:operateur is gateau diff}, we easily obtain Lemma \ref{lemma:existence of martingale coupling}.

\smallskip

\noindent\textbf{Proof of Lemma \ref{lemma:existence of martingale coupling}.} 
We verify that the mapping $ \Uc$ satisfies Assumption \ref{ass:gen ext implicit function theorem}. 
Consider $\mathbb{L}^{\infty}_0 (\mu_1) $ equipped with the norm $\Vert \cdot \Vert_{\mathbb{L}^{\infty} (\mu_1)}$ and $ \mathbb{L}^{\alpha}_0 (\mu_1) $ endowed with the norm $\Vert \cdot \Vert_{\mathbb{L}^{\alpha} (\mu_1)} $. 
By Hölder's inequality, $\Vert \cdot \Vert_{\mathbb{L}^{\alpha} (\mu_1)} \leq \Vert \cdot \Vert_{\mathbb{L}^{\infty} (\mu_1)} $ and, by a standard argument in measure theory, $ a \in \mathbb{L}^\alpha_0 (\mu_1) \mapsto \Vert a \Vert_{\mathbb{L}^{\infty} (\mu_1)} $ is lower semi-continuous. 
Now notice that $ \mathrm{d}_a \Uc (0, 0) = {\rm Id} - \E^{\mu}_1 \circ \E^{\mu}_2 $.
Thus, by Lemmas \ref{lemma:closedness ans coercivity}, \ref{lemma:operateur is gateau diff} and \ref{lemma:continuity gateau diff}, $\Uc$ satisfies Assumption \ref{ass:gen ext implicit function theorem}, proving the desired result. \ep

\medskip

\noindent \textbf{Proof of Lemma \ref{lemma:closedness ans coercivity}.} We will prove later that $ \vertiii{\E^{\mu}_1 \circ \E^{\mu}_2 }_\alpha < 1$. Assume for now that it is true.
Again, we will extensively use the identification $ \mathbb{L}^\alpha (\mu_i) \subset \mathbb{L}^\alpha (\mu) $.

\noindent{\bf Proof of point} \ref{direct_sum_sumspace}. Let $ (f_1 , f_2) \in \mathbb{L}^\alpha_0 (\mu_1) \times \mathbb{L}^\alpha (\mu_2) $ be such that $ f_1 \oplus f_2 = 0 $. 
Taking conditional expectations with respect to $X_1$ and with respect to $X_2$, we get the two following equations
$$
f_1 (X_1) + \E^{\mu}_1 [f_2 ] = 0 \,\,\,, \,\,\, f_2 (X_2) + \E^{\mu}_2 [f_1 ] = 0
.$$
Substituting $f_2 = - \E^{\mu}_2 [f_1 ]$ in the first equation, we get $ f_1(X_1) - \E^{\mu}_1 \circ \E^{\mu}_2 (f_1) = 0 $. 
Since $ \vertiii{ \E^{\mu}_1 \circ \E^{\mu}_2}_\alpha < 1$, $ {\rm Id} - \E^{\mu}_1 \circ \E^{\mu}_2$ is invertible, so $f_1 = 0$. 
Hence, as $f_1 \oplus f_2 = 0$, we also have $f_2 = 0$, proving that the sum is direct.

\noindent We now move on to the closure property. 
Let $ u_n:= f_1^n \oplus f_2^n $ be such that $u_n \rightarrow u $ in $\mathbb{L}^\alpha (\mu) $. 
Then, by similar considerations, we obtain $ f_1^n -\E^{\mu}_1 \circ \E^{\mu}_2 (f_1^n) = \E^{\mu}_1 [ u_n -\E^{\mu}_2 [u_n] ] $. 
Hence, since $ {\rm Id} - \E^{\mu}_1 \circ \E^{\mu}_2$ is invertible, and the conditional expectation is continuous, $f^n_1 \rightarrow ({\rm Id} - \E^{\mu}_1 \circ \E^{\mu}_2)^{-1} (\E^{\mu}_1 [ u -\E^{\mu}_2 [u] ]) $ which also proves that $f_2^n $ converges to some $f_2$ and that $ u = f_1 \oplus f_2 $. 

\noindent{\bf Proof of point} \ref{coercivity direct sum}.
We proved that $ \mathbb{L}^\alpha_0 (\mu_1) \oplus \mathbb{L}^\alpha (\mu_2) $ is direct and closed. 
Hence, there exists $C > 0$ such that for all $ f_1 \oplus f_2 \in \mathbb{L}^\alpha_0 (\mu_1) \oplus \mathbb{L}^\alpha (\mu_2) $, we have 
$$
\Vert f_1 \oplus f_2 \Vert_{\mathbb{L}^{\alpha } (\mu) } \geq C (\Vert f_1 \Vert_{\mathbb{L}^{\alpha } (\mu) } + \Vert f_2 \Vert_{\mathbb{L}^{\alpha } (\mu) }) 
.$$

\smallskip

\noindent{\bf Proof of point} \ref{contraction_operator}. Let $ 1 \leq \alpha \leq \infty $. 
We will now prove that $ \E^{\mu}_1 \circ \E^{\mu}_2: \mathbb{L}^\alpha_0 (\mu_1) \rightarrow \mathbb{L}^\alpha_0(\mu_1) $ is a contraction. 
We distinguish two cases.

\noindent{\bf Case $1$}: $\alpha < \infty $.
Assume, to the contrary that $\vertiii{\E^{\mu}_1 \circ \E^{\mu}_2}_1 = 1$. 
We first prove that there exists $ u \in \mathbb{L}^{1}_0 (\mu_2) $ such that $ \E^\mu [ \vert u (X_2) \vert ] \leq 1 $ and $ \E^{\mu_1} \big[ \vert \E^{\mu}_1 [ u ] \vert \big] = 1$.
Let $(v_n)_n \in \mathbb{L}^{1}_0 (\mu_1) $ be a sequence of functions such that $ \Vert v_n \Vert_{\mathbb{L}^{\alpha} (\mu)} \leq 1 $ and $ \Vert \E^{\mu}_1 \big[ \E_2^\mu \big[ v_n]] \Vert_{\mathbb{L}^{\alpha} (\mu)} \rightarrow 1$.
Since $q_2 (X_2) \geq c $, and $\kappa$ is bounded by $C$, a quick computation yields 
$$
\vert u_n (X_2) \vert := \vert \E^{\mu}_2 [ v_n (X_1) ] \vert  = \int_\R v_n(x_1) \frac{\kappa(x_1, X_2)}{q_2 (X_2)} \mu_1 (\mathrm{d}x_1) \leq \frac{C}{c} \int_\R \vert v_n(x_1) \vert \mu_1 (\mathrm{d}x_1) \leq  \frac{C}{c}
.$$
This proves that $u_n \in \mathbb{L}^\infty (\mu_2) $, with $ \Vert u_n \Vert_{\mathbb{L}^\infty (\mu_2)} \leq C/c $ for all $n \in \N$. 
Hence $(u_n)$ is also bounded in $\mathbb{L}^2 (\mu_2)$. 
By weak compactness, there exists $ u \in \mathbb{L}^{2}(\mu_2) $ such that 
$
u_n \rightharpoonup u .$ 
In other terms, for all $v \in \mathbb{L}^{2}(\mu_2) $, we have 
\begin{equation}\label{eq:weak cone}
\int_\R u_n (x_2) v(x_2) q_2 (x_2) \mathrm{d}x_2 \rightarrow \int_\R u (x_2) v(x_2) q_2 (x_2) \mathrm{d}x_2 ,
\end{equation}
or equivalently $
\int_\R (u_n  - u ) (x_2) v(x_2) q_2 (x_2) \mathrm{d}x_2 \rightarrow 0
.$
Since $q_2 (x_2) \geq c $ on the support of $\mu_2$, and $ \kappa $ is essentially bounded, we have $ \int_\R \int_\R \Big\vert \frac{\kappa (x_1, x_2)}{ q_2 (x_2)} \Big\vert^2 q_2(x_2) \mathrm{d}x_2 \mu_1(\mathrm{d}x_1) <\infty $. 
Therefore, by Fubini's theorem, $ \mu_1 $-almost surely, $x_2 \mapsto \frac{\kappa (x_1, x_2)}{ q_2 (x_2)} $ is in $\mathbb{L}^{2}(\mu_2) $, hence, by equation \eqref{eq:weak cone}, we have 
$$
\int_\R (u_n (x_2) - u (x_2)) \kappa (X_1, x_2) \mathrm{d}x_2 \rightarrow 0
\,\, \mu_1-\textit{a.s.}
.$$
Also, since $ \Vert u_n \Vert_{\mathbb{L}^\infty (\mu) } \leq C $, the sequence is $\alpha$-uniformly integrable; thus 
$$
\E_1^\mu [ u_n ] = \int_\R u_n (x_2) \kappa (X_1, x_2) \mathrm{d}x_2 \xrightarrow[n \rightarrow \infty]{\mathbb{L}^{\alpha} (\mu_1)} \E_1^\mu [ u ] = \int_\R u (x_2) \kappa (X_1, x_2) \mathrm{d}x_2
.$$
Since $ \Vert \E^{\mu}_1 [ u_n ] \Vert_{\mathbb{L}^{\alpha } (\mu) }= \Vert \E^{\mu}_1[ \E^{\mu}_2 [ v_n ]] \Vert_{\mathbb{L}^{\alpha} (\mu) } \rightarrow 1 $, this proves that $ \Vert \E^{\mu}_1 [ u ] \Vert_{\mathbb{L}^{\alpha } (\mu) } = 1$. 
Furthermore, since the conditional expectation is a contraction, 
$\Vert u_n \Vert_{\mathbb{L}^{\alpha } (\mu) } \leq \Vert v_n \Vert_{\mathbb{L}^{1 } (\mu) } \leq 1 $ hence, $ \Vert u \Vert_{\mathbb{L}^{\alpha } (\mu) } \leq 1$. 
Finally, notice that $\E^{\mu}[ u_n ] = 0 $, which implies that $ \E^{\mu} [ u ] = 0 $ as $n$ goes to infinity. 
In this case, $\vert\E^{\mu}_1 [u ] \vert 
=
\E^{\mu}_1 [ \vert u \vert ]
$, which, by the equality case of Jensen's inequality for $x \mapsto \vert x \vert^\alpha$, implies that $ \text{sgn} (u) $ is $\sigma(X_1)$-measurable. 
Since $\sigma(X_1) \cap \sigma(X_2)$ is trivial by Assumption \hyperref[ass:mart coupling]{ $ \textbf{\rm C}_{\W^{\rm ad}_p}$} \ref{cond:info discrepancy}, $u$ must have a constant sign, which is $0$ since $\E^{\mu} [ u ] = 0$—a contradiction with $\Vert u \Vert = 1$.

\noindent{\bf Case $2$}: $\alpha = \infty$.
Assume, to the contrary, that $ \vertiii{\E^{\mu}_1 \circ \E^{\mu}_2}_\infty = 1 $. 
Then there exists a sequence $ (u_n)_n \subset \mathbb{L}^\infty_0(\mu_1) $ such that 
$$
\Vert u_n \Vert_{\mathbb{L}^\infty(\mu_1)} \leq 1 \quad \text{and} \quad \Vert \E^{\mu}_1[\E^{\mu}_2[u_n]] \Vert_{\mathbb{L}^\infty(\mu_1)} \longrightarrow 1.
$$
Since $(u_n)$ is uniformly bounded in $\mathbb{L}^\infty(\mu_1) $, it is also bounded in $\mathbb{L}^2(\mu_1)$, and thus (up to a subsequence) $u_n \rightharpoonup u $ weakly in $\mathbb{L}^2(\mu_1) $, for some $u \in \mathbb{L}^2(\mu_1) $. 
Using the disintegration of $ \mu $ and similar computations as in the case $ \alpha = 1 $, we obtain
$$
\E_2^\mu[u_n] \longrightarrow \E_2^\mu[u] \quad \text{in } \mathbb{L}^2(\mu_2).
$$
Moreover, since $ \kappa $ is bounded and $ q_2 $ is strictly positive on the support of $ \mu_2 $, we derive the key estimate:
\begin{equation}\label{ineq:cond-expectation-estimate}
\left| \E_1^\mu[\E_2^\mu[u_n]] - \E_1^\mu[\E_2^\mu[u]] \right| \leq C \left\| \E_2^\mu[u_n] - \E_2^\mu[u] \right\|_{\mathbb{L}^1(\mu_2)}.
\end{equation}
Hence, we conclude:
$$
\left\| \E_1^\mu[\E_2^\mu[u_n]] - \E_1^\mu[\E_2^\mu[u]] \right\|_{\mathbb{L}^\infty(\mu_1)} \longrightarrow 0,
$$
and in particular, $ \Vert \E_1^\mu[\E_2^\mu[u]] \Vert_{\mathbb{L}^\infty(\mu_1)} = 1 $. 
Since the unit ball of $ \mathbb{L}^\infty(\mu_1) $ is convex and closed in $ \mathbb{L}^2(\mu_1) $, it is weakly closed. 
Thus $ \Vert u \Vert_{\mathbb{L}^\infty(\mu_1)} \leq 1 $, and in fact equality must hold, $ \Vert u \Vert_{\mathbb{L}^\infty(\mu_1)} = 1 $, due to the contraction property of conditional expectations.
Define
$$
Y_n:= |u|^n \in \sigma(X_1), \qquad Z_n:= |\E_2^\mu[u]|^n \in \sigma(X_2).
$$
Since $ \Vert u \Vert_{\mathbb{L}^\infty(\mu_1)} = 1 $ and $ \Vert \E_2^\mu[u] \Vert_{\mathbb{L}^\infty(\mu_2)} = 1 $, we have:
$$
Y_n \xrightarrow[n \to \infty]{\textit{a.s.}} \mathds{1}_{\{|u| = 1\}}, \qquad Z_n \xrightarrow[n \to \infty]{\textit{a.s.}} \mathds{1}_{\{|\E_2^\mu[u]| = 1\}}.
$$
Since both $ Y_n $ and $ Z_n $ are bounded in $ \mathbb{L}^\infty $, the convergence also holds in $ \mathbb{L}^1(\mu) $. 
Combining this with estimate~\eqref{ineq:cond-expectation-estimate}, we obtain:
$$
\Vert \E_1^\mu[\mathds{1}_{\{|\E_2^\mu[u]| = 1\}}] \Vert_{\mathbb{L}^\infty(\mu_1)} = \Vert \E_1^\mu[\E_2^\mu[\mathds{1}_{\{|u| = 1\}}]] \Vert_{\mathbb{L}^\infty(\mu_1)} = \Vert \E_2^\mu[\mathds{1}_{\{|u| = 1\}}] \Vert_{\mathbb{L}^\infty(\mu_2)} = 1.
$$
Equivalently,
$$
\text{ess\,inf } \E_1^\mu[\mathds{1}_{\{|\E_2^\mu[u]| < 1\}}] = \text{ess\,inf } \E_1^\mu[\E_2^\mu[\mathds{1}_{\{|u| < 1\}}]] = \text{ess\,inf } \E_2^\mu[\mathds{1}_{\{|u| < 1\}}] = 0.
$$
By Assumption \hyperref[ass:mart coupling]{ $ \textbf{\rm C}_{\W^{\rm ad}_p}$} \ref{cond:disintegration}, as $ \kappa $ is greater than some $c >0$ $\mu$ almost surely, this implies that:
$$
\mathds{1}_{\{|u| < 1\}} = 0 \quad \text{and} \quad \mathds{1}_{\{|\E_2^\mu[u]| < 1\}} = 0, \quad \mu -\textit{a.s.}
$$
Thus, $ |u| = 1 $ and $ |\E_2^\mu[u]| = 1 $ $ \mu $-almost surely. 
As in the case $\alpha=1$, this is only possible if $ u $ is constant (equal to $ 1 $ or $ -1 $), which contradicts $\E^\mu[u] = 0 $. 
Hence
$$
\vertiii{\E^{\mu}_1 \circ \E^{\mu}_2}_\infty < 1.
$$
\ep

\subsection{An Implicit Function Theorem} \label{subsec:implicit function theorem}

For the purpose of our proofs, we need to adapt the Implicit Function Theorem in order to dispense with the usual $C^1$ regularity, which is lost in our setting. 
This situation is reminiscent of the well-known case of Nemytskij operators. 
A variety of extensions of the Implicit Function Theorem can be found in the literature: some consider non-Banach spaces, as in \citeauthor{hamilton1982inverse} \cite{hamilton1982inverse} or \citeauthor{ekeland2011inverse} \cite{ekeland2011inverse}; others weaken the regularity assumptions on the functional, as in \citeauthor{biasi2008implicit} \cite{biasi2008implicit} and \citeauthor{accinelli2009generalization} \cite{accinelli2009generalization}.
However, none of these results apply directly to our framework, since we wish to apply the theorem to the functional $\mathcal{U}$ defined by \eqref{eqdef:operateur}, for which we only have the properties established in Lemmas \ref{lemma:closedness ans coercivity}, \ref{lemma:continuity gateau diff}, and \ref{lemma:operateur is gateau diff}. 
Our approach is close to \citeauthor{wachsmuth2014differentiability} \cite{wachsmuth2014differentiability} where others proved a weaker form of the implicit function theorem by considering inclusion of spaces. However, they assume the existence of an implicit mapping, which we want to prove here.

Let $ (E_1, \Vert \cdot \Vert_1) $ and $ (E_2, \Vert \cdot \Vert_2)$ be two Banach spaces, with $ E_1 \subset E_2$, $ I $ an open subset of $\R$ and $(t_0, a_0) \in I \times E_1 $. 
Let $ \Psi: I \times E_2 \rightarrow E_1 $ be a mapping, using notations $\Lc(E_i)$ and $\vertiii{ \cdot }$ defined in Equation \eqref{eqdef:continuousoperator} and \eqref{eqdef:def operator norm}, we define the following set of Assumptions.

\begin{Assumption}\label{ass:gen ext implicit function theorem}
$E_1$, $E_2$ and $ \Psi$ satisfy:
\begin{enumerate}[label=\textnormal{(\roman*)}]
\item \label{norm inequality} $\Vert \cdot \Vert_2 \leq \Vert \cdot \Vert_1$ and $ \Vert \cdot \Vert_1: E_2 \rightarrow \bar{\R}^+$ is lower semi-continuous.
\item \label{Gâteau differentiability} $ \Psi$ is Gâteaux differentiable in a neighborhood of $(t_0, a_0)$ in $I \times E_2 $, with Gâteaux derivative with respect to $a$, $ \mathrm{d}_a \Psi (t, a) \in \Lc_c (E_2)\cap \Lc_c (E_1)$.
The following mapping is continuous 
$$ \lambda \in [0, 1] \mapsto \mathrm{d}_a \Psi (t, a' + \lambda (a - a')) (a - a') \in E_2,$$ 
and $ t \mapsto \partial_t \Psi (t, a) \in E_2$ is continuous in a neighborhood of $ (t_0, a_0)$.
\item \label{Gâteau diff at 0} The operator $ S:= \mathrm{d}_a \Psi (t_0, a_0) - { \rm Id}$ satisfies $ \vertiii{ S }_2, \vertiii{ S }_1 < 1$. 
\item \label{continuity at 0 Gâteau diff} For $i =1,2$, the mappings $\mathrm{d}_a \Psi: I \times E_1\rightarrow  (\Lc_c (E_i), \vertiii{\cdot} )$ are continuous at $ (t_0, a_0)$.
\end{enumerate}
\end{Assumption}

\begin{Lemma}\label{lemma:revised implicit function theorem}
Under Assumption \ref{ass:gen ext implicit function theorem}, there exists $ \eta > 0$, and $ a: (- \eta, \eta) \longrightarrow V_{ \eta}:= \{ a \in E_1 \, \, : \,\, \Vert a \Vert_1 \leq \eta \} $ such that for all $ \vert t - t_0 \vert < \eta$, we have 
\begin{equation*}
 \Psi (t, a(t)) = \Psi(t_0, a(t_0)) \,\, \text{and} \,\,
\Vert a(t) -  a(t_0) \Vert_2 \leq C \vert t-t_0 \vert (\Vert \partial_t \Psi (t_0, a_0) \Vert_2 + o (1) ) 
.\end{equation*}
\end{Lemma}

\noindent\textbf{Proof of Lemma \ref{lemma:revised implicit function theorem}.} 
This proof is an adaptation of the Implicit Function Theorem. 
Without loss of generality, assume that $ a_0 = 0 $, $ t_0 = 0$ and $\Psi (0, 0) = 0$.
Define $ \Gc (t, a):= a - \Psi (t, a) $. 
Let $\varepsilon > 0$. 
By the continuity of $\mathrm{d}_a \Psi$ from Assumption \ref{ass:gen ext implicit function theorem} \ref{continuity at 0 Gâteau diff}, there exists $ \eta > 0$ (the dependence on $\varepsilon $ is omitted for the sake of clarity) such that for all $ \vert t \vert < \eta $ and $ \Vert a \Vert_1 \leq \eta $, we have $\vertiii{ \mathrm{d}_a \Psi (t, a) - \mathrm{d}_a \Psi (0, 0) }_1 \leq \varepsilon $ and $\vertiii{ \mathrm{d}_a \Psi (t, a) - \mathrm{d}_a \Psi (0, 0)}_2 \leq \varepsilon $.
Let $ V_\eta$ be defined as in the statement of the lemma.
Then $V_\eta $ is convex and, by lower semi-continuous Assumption \ref{ass:gen ext implicit function theorem} \ref{norm inequality}, $V_\eta $ is closed in $E_2$. 
Let $a_1, a_2 \in V_\eta$ and $ \vert t \vert < \eta$, and set $a^\lambda := a_2 + \lambda (a_1 - a_2)$ for $\lambda \in [0, 1]$. 
By Assumption \ref{ass:gen ext implicit function theorem} \ref{Gâteau differentiability}, $ \Psi $ is continuously differentiable on the segment $ [a_1, a_2] $, and by Assumption \ref{ass:gen ext implicit function theorem} \ref{Gâteau diff at 0}
\begin{align*}
 \Gc (t, a_1) - \Gc(t, a_2) &= (a_1 - a_2) - \int_0^1 \mathrm{d}_a \Psi (t, a^\lambda ) (a_1 - a_2) \mathrm{d}\lambda
 \\
 &=
 -S (a_1 - a_2) -
 \int_0^1 
 \Big(\mathrm{d}_a \Psi (t, a^\lambda ) - \mathrm{d}_a \Psi (0, 0) \Big) (a_1 - a_2) \mathrm{d}\lambda
.\end{align*}
Now, since $ V_\eta$ is convex, we have that for all $ 0\leq \lambda \leq 1$, $a^\lambda \in V_\eta$.
Furthermore, since $ \Vert \cdot \Vert_1$ is convex and lower semi-continuous by Assumption \ref{ass:gen ext implicit function theorem} \ref{norm inequality}, we have the triangle inequality
$$
 \big\Vert 
 \int_0^1 
 \Big(\mathrm{d}_a \Psi (t, a^\lambda ) - \mathrm{d}_a \Psi (0, 0) \Big) (a_1 - a_2) \mathrm{d}\lambda
\big\Vert_i 
\leq 
\varepsilon \Vert a_1 - a_2 \Vert_i.
$$
Furthermore, $\lambda \mapsto \big(\mathrm{d}_a \Psi (t, a^\lambda ) - \mathrm{d}_a \Psi (0, 0) \big) (a_1 - a_2) $ is continuous by Assumption \ref{ass:gen ext implicit function theorem} \ref{Gâteau differentiability}; we get, for $a_2 =0$, 
\begin{equation}\label{ineq:stable set 1}
\Vert \Gc (t, a_1) \Vert_1 
\leq 
\vertiii{S}_1 \Vert a_1 - a_2 \Vert_1 + \varepsilon \Vert a_1 - a_2 \Vert_1
=
(\varepsilon + \vertiii{S}_1) \Vert a_1 \Vert_1
.\end{equation}
By similar considerations, we obtain 
\begin{equation}\label{ineq:stable set 2}
\Vert \Gc (t, a_1) - \Gc(t, a_2) \Vert_2
\leq 
\vertiii{S}_2 \Vert a_1 - a_2 \Vert_2 + \varepsilon \Vert a_1 - a_2 \Vert_2
=
(\varepsilon + \vertiii{S}_2) \Vert a_1 - a_2 \Vert_2
.\end{equation}
Both inequalities \eqref{ineq:stable set 1} and \eqref{ineq:stable set 2} ensure that, for $ \varepsilon > 0 $ small enough that $(\varepsilon + \vertiii{S}_i) < 1 $ for $ i= 1, 2$ (which is possible since $ \vertiii{S}_i < 1 $ by Assumption \ref{ass:gen ext implicit function theorem} \ref{Gâteau diff at 0}), and $\vert t \vert < \eta$, we have $ \Gc: V_\eta \rightarrow V_\eta$ is well defined and a contraction.
Hence, by the Banach--Picard fixed-point theorem, since $ V_\eta$ is closed in a Banach space, there exists $ a: (- \eta, \eta) \rightarrow V_{ \eta}$ such that for all $ \vert t \vert < \eta$, $\Psi (t, a(t)) = \Psi(0, 0) $. 

Now, by Assumption \ref{ass:gen ext implicit function theorem} \ref{Gâteau differentiability}, $\Psi$ is continuously differentiable with respect to $t$, so we have:
 $$
 \int_0^1 \mathrm{d}_a \Psi (t, \lambda a(t)) (a(t)) \mathrm{d} \lambda = -t \int_0^1 \partial_t \Psi (\lambda t, 0) \mathrm{d} \lambda
 .$$
 This can be rewritten as 
 $$
 ({ \rm Id} + S) \frac{a(t)}{t} = -\int_0^1 \Big(\mathrm{d}_a \Psi (t, \lambda a(t)) - \mathrm{d}_a \Psi (0,0) \Big) \Big(\frac{a(t)}{t} \Big) \mathrm{d} \lambda - \int_0^1 \partial_t \Psi (\lambda t, 0) \mathrm{d} \lambda
. $$
Now, since $ \vertiii{S}_2 < 1$, $ { \rm Id} + S$ is invertible. Furthermore, as $a(t) \in V_\eta$, we have $ \Big\Vert \big(\mathrm{d}_a \Psi (t, \lambda a(t)) - \mathrm{d}_a \Psi (0,0) \big) \Big(\frac{a(t)}{t} \Big) \Big\Vert_2 \leq \varepsilon \Big\Vert \frac{a(t)}{t} \Big\Vert_2 $ hence 
$$
\Big\Vert \frac{a(t)}{t} \Big\Vert_2
\leq 
\vertiii { (\rm Id + S)^{-1}}_{2} \big(\varepsilon \Big\Vert \frac{a(t)}{t} \Big\Vert_2 + \Vert \partial_t \Psi ( 0, 0) \Vert_2 + o (1) \big) 
.$$
This provides for sufficiently small $\varepsilon$
$$
\Big\Vert \frac{a(t)}{t} \Big\Vert_2 
\leq 
\frac{1}{ 1 - \varepsilon \vertiii{ (\rm Id + S)^{-1} }_{2}} \big( \Vert \partial_t \Psi ( 0, 0) \Vert_2 + o (1) \big) ,
$$
which proves the desired inequality, since $\frac{1}{ 1 - \varepsilon \vertiii{ (\rm Id + S)^{-1} }_{2} }$ can be bounded by a constant independently of $\varepsilon$. 
\ep 

\subsection{Expansions in Probability and $\mathbb{L}^p(\mu)$}\label{Expansions in Probability and}

\begin{Lemma}\label{lemma:expans in Lp C1}
Let $ (Z_r)_{ r > 0} $ be a family of $I$-valued random variables, for some open interval $I$, on $(\Omega, \Fc, \P)$.
Assume that there exist two random variables $Z_0$ and $U$ such that
$$
Z_r = Z_0 + rU + o_{\mathbb{L}^p (\P)}(r)
.$$
Assume that there exists $ A \subset \Omega $ and $K \subset \R$ an interval, such that $ Z_r = Z_0 $ on $A$, $U = 0$ on $A$ and $ Z_r, Z_0 \in K $ on $A^{\rm c}$.
Let $ F: I \rightarrow \R $ be in $ C^1 (I, \R) $, such that $ \Vert F'\Vert_{\mathbb{L}^{\infty}(K)} < \infty $, then 
$$
F(Z_r) = F(Z_0) + rF'(Z_0) U + o_{\mathbb{L}^p (\P)}(r).
$$
\end{Lemma}

\proof  
Since $I$ is an open interval and $ \P (Z_r \in I) = 1$ for all $r \geq 0$, we can apply Taylor's formula,
$$
F(Z_r) = F(Z_0) + (Z_r - Z_0)F'(Z_0) + (Z_r - Z_0) \int_0^1  \big( F' (Z_0 + \lambda (Z_r - Z_0)) -  F' (Z_0 ) \big) \mathrm{d} \lambda.
$$
Letting $\Delta_r :=  F(Z_r) - F(Z_0) - r F'(Z_0) U $, we have 
\begin{align*}
\vert \Delta_r \vert 
\leq
\vert F' (Z_0) \vert \vert Z_r - Z_0 - rU \vert + 
\vert Z_r - Z_0\vert \int_0^1 \vert F' (Z_0 + \lambda (Z_r - Z_0) ) - F'(Z_0) \vert \mathrm{d}\lambda
.
\end{align*}
Since $F'$ is bounded on $K$, $ \vert F'(Z_0) \vert \mathds{1}_{A^{\rm c}} \leq C := \text{sup}_{ x \in K } \vert F'(x) \vert$.
Furthermore $\vert Z_r - Z_0 - rU \vert  \mathds{1}_{A} = 0 $ hence 
\begin{equation}\label{op1 terme 1}
    \vert F' (Z_0) \vert \vert Z_r - Z_0 - rU \vert = o_{\P} (r).
\end{equation}
Notice that
\begin{align*}
&\vert Z_r - Z_0\vert \int_0^1 \big\vert F'\big(Z_0 + \lambda (Z_r - Z_0) \big) - F'(Z_0) \big\vert \mathrm{d}\lambda
\\
&=
\mathds{1}_{A^{\rm c}}
 \vert Z_r - Z_0\vert \int_0^1 \big\vert F'\big(Z_0 + \lambda (Z_r - Z_0) \big) - F'(Z_0) \big\vert \mathrm{d}\lambda ,
\end{align*}
since $ Z_r= Z_0  $ on $A$.
Furthermore, since $ Z_r, Z_0 \in K $ on $A^{\rm c}$ and $K$ is an interval, for all $\lambda \in [0, 1]$, $Z_0 + \lambda (Z_r - Z_0) \in K$ hence 
$$\mathds{1}_{A^{\rm c}}
\vert F' (Z_0 + \lambda (Z_r - Z_0) ) - F'(Z_0) \vert
\leq 2 \Vert F'\Vert_{\mathbb{L}^{\infty}(K)}\mathds{1}_{A^{\rm c}}
.$$
Hence, by \eqref{op1 terme 1}, we get 
\begin{align}\label{eqdef:estimatedeltarconvergenceproba}
\vert\Delta_r \vert
\leq 
2 r \vert U \vert  \Vert F'\Vert_{\mathbb{L}^{\infty}(K)}\mathds{1}_{A^{\rm c}}
+ 
2\vert Z_r - Z_0 - rU \vert  \int_0^1 \vert F' ( \bar{Z}^\lambda_r ) - F'(Z_0) \vert \mathrm{d}\lambda
\mathds{1}_{A^{\rm c}}
+
o_{\P}(r),
\end{align}
where $ \bar{Z}^\lambda_r := Z_0 + \lambda (Z_r - Z_0)$.
Now, notice that $(Z_r, Z_0) \xrightarrow[ r \rightarrow 0]{\P} (Z_0, Z_0) $, and the map 
$H : (x, y ) \in K \times K \mapsto \int_0^1 \vert F' ( y + \lambda(x - y) ) - F'(y) \vert \mathrm{d}\lambda $ is continuous since $F'$ is continuous and bounded on $K$, with $H(y, y) = 0$.
So, by the continuous mapping theorem, we have
$$
\int_0^1 \vert F' ( \bar{Z}^\lambda_r ) - F'(Z_0) \vert \mathrm{d}\lambda
= H (Z_r, Z_0)
\xrightarrow[ r \rightarrow 0]{\P} H(Z_0, Z_0) = 0
.$$
Now putting the last convergence result along with \eqref{eqdef:estimatedeltarconvergenceproba} yields
$$ F(Z_r) = F(Z_0) + rF'(Z_0) U + o_{\P}(r)
.$$
To obtain the $\mathbb{L}^p$ estimate, it suffices to prove that the family $ ( \vert \Delta_r / r \vert^p)_r $ is uniformly integrable. 
This follows from $ \Delta_r \mathds{1}_A = 0 $ and from the following bound: setting $C := \text{sup}_{ x \in K } \vert F'(x) \vert$, since $\P (Z_r  \mathds{1}_{A^{\rm c}} \in K , Z_0  \mathds{1}_{A^{\rm c}} \in K ) = 1$, we have 
$$
\vert \Delta_r \vert \mathds{1}_{A^{\rm c}}  \leq C \vert Z_r - Z_0 \vert \mathds{1}_{A^{\rm c}} +
 r\vert U\vert \Vert F' \Vert_{ \mathbb{L}^\infty (K) }.
$$
Now, since $\frac{Z_r - Z_0}{r} \rightarrow U \in \mathbb{L}^p (\P)$, $(\frac{Z_r - Z_0}{r})_{r}$ is uniformly integrable and we obtain the desired result. 
\ep

\medskip

\noindent In the following, $(\Zc, \Bc, \nu)$ denotes a finite-dimensional normed vector space endowed with its Borel $\sigma$-algebra and a probability measure $\nu$, $ 1 \leq \alpha \leq +\infty$, and $(\Theta_n)_n$ is a sequence of functions $\Theta_n : \R \times \Zc \rightarrow \R$.
\begin{Assumption}\label{ass:change of variable}
\begin{enumerate}[label=\textnormal{(\roman*)}]
\item \label{cond:measurability in z}  For all $x \in \R$, the mapping $ \Theta_n( x, \cdot )$ is measurable. 
\item \label{cond:regularity in x} For $\nu$ almost every $ z \in \Zc$, the mapping $ \Theta_n ( \cdot, z) $ is $C^2$.
\item  \label{cond:convergence} There exist $C>0$, a sequence of positive real numbers, $(C_n)_n $ and a sequence of functions $(f_n)_n$, with $\vert \Theta_n (x, z ) \vert \leq f_n (z)$, $\vert \partial_x \Theta_n (x, z ) \vert \leq C_n$ and $\vert \partial_{xx}^2 \Theta_n (x, z) \vert \leq C$,  where $ f_n \in \mathbb{L}^{\alpha}(\nu)$. Furthermore, the sequences satisfy $ C_n \rightarrow 0 $ and $ \Vert f_n \Vert_{\mathbb{L}^{\alpha} (\nu )} \rightarrow 0$.
\end{enumerate}
\end{Assumption}

\begin{Lemma}\label{lemma:estim diffeo and change of variable}
Let $(\Zc, \Bc, \nu)$, $\alpha$ and $(\Theta_n)_n$ be as above, such that Assumption \ref{ass:change of variable} holds.
Then there exists $N \in \N$ and a family of functions $ ( \varphi_n )_{n \geq N}$, $ \varphi_n : \R \times \Zc \rightarrow \R $ such that 
\begin{equation}\label{eqdef:phi change of variable}
x + \Theta_n (x, z) = y
\,\,
\text{iff} \,\,
x  = \varphi_n (y, z) 
\,\,
\text{ for all } n \geq N \,\, , \,\, x,y  \in \R \,\, \text{and} \,\, z \in \Zc
.\end{equation}
Furthermore, $\varphi_n$ satisfies the following 
\begin{enumerate}[label=\textnormal{(\roman*)}]
\item  For all $y \in \R$, and $n \in \N$, the mapping $ \varphi_n (y, \cdot)$ is measurable. 
\item  For $\nu$-almost $ z \in \Zc$, the mapping $ \varphi_n ( \cdot, z) $ is $C^2$. There exists $C>0$ such that for all $y \in \R$ and $\nu$ almost-every $z$, $\vert \partial_y \varphi_n (y, z) \vert \leq C$ and $\vert \partial_{yy}^2 \varphi_n (y, z) \vert \leq C$. Furthermore, 
\begin{align*}
\partial_{y} \varphi_n (y, z) = \frac{ 1}{1 + \partial_{x} \Theta_n (\varphi_n (y, z), z)}
.\end{align*}
\item The following convergence holds
\begin{align*}
\text{max} ( u_n, v_n )  \xrightarrow [ n \rightarrow +\infty] {} 0 ,
\end{align*}
where $ u_n :=\Vert {\rm sup}_{y \in \R} \{ \vert \varphi_n (y, \cdot) - y \vert \} \Vert_{\mathbb{L}^{\alpha}(\nu)}$, $v_n :=
\Vert {\rm sup}_{y \in \R} \{ \vert \partial_y \varphi_n (y, \cdot) - 1 \vert \} \Vert_{\mathbb{L}^{\alpha}(\nu)}
$.
\end{enumerate}
\end{Lemma}

\noindent \textbf{Proof.} 
Fix $n \in \N $, up to a multiplication by an indicator, we can assume in the rest that $\Theta_n$ is $C^1$ in $x$ for every $z \in \Zc$. 
Define the following sequence of functions 
$
u_{k+1, n}(y, z) := y - \Theta_n(u_{k, n} (y, z) , z)
$ with $u_{0, n}(y, z) = y$. 
By Condition \ref{cond:convergence} of Assumption \ref{ass:change of variable}, $ \vert \partial_x \Theta_n( x, z) \vert \leq C_n \rightarrow 0$.
Hence there exists $ N \in \N$ such that, for $n \geq N$, there is $ 0 \leq \eta < 1$ such that 
\begin{equation}
\vert u_{k+1, n}- u_{k, n} \vert (y, z)  \leq \eta \vert u_{k, n}- u_{k-1, n} \vert (y, z).
\end{equation}
Thus there exists $\varphi_n :\R \times \Zc \rightarrow \R  $, such that $u_{k, n} (y, z) \xrightarrow[k \rightarrow \infty]{} \varphi_n (y, z)$, proving measurability of $\varphi_n$ as a pointwise limit of measurable functions.
Let $v_{k, n} := \partial_y u_{k, n} (y, z)$. Then 
\begin{align*}
v_{k+1, n} (y, z) = 1 &- \partial_{x}\Theta_n  (\varphi_n (y, z), z ) v_{k, n} (y, z) 
\\
&+ \big( \partial_{x}\Theta_n (\varphi_n (y, z), z ) - \partial_{x}\Theta_n  (u_{k, n} (y, z) , z ) \big) v_{k, n} (y, z)
.\end{align*}
By Condition \ref{cond:convergence} of Assumption \ref{ass:change of variable}, $ \vert \partial_x \Theta_n( x, z) \vert \leq C_n \rightarrow 0$ so $\vert v_{k} (y, z) \vert \leq C_1$, for some $C_1 >0$. 
Furthermore, since $\vert \partial_{xx}^2 \Theta_n \vert \leq C_2$, for some $C_2$, we obtain, for $n$ large enough (uniformly in $y, z$),  
$$
\vert v_{k+1}- v_{k} \vert (y, z)  \leq \eta \vert v_{k}- v_{k-1} \vert (y, z) + C_1C_2 ( \vert u_{k}- \varphi_n \vert (y, z) +\vert u_{k-1}- \varphi_n \vert (y, z)   ) 
.$$
And so, $ v_k (y, z ) \xrightarrow[k \rightarrow \infty]{}\partial_y \varphi_n (y,z) = \frac{1}{1 + \partial_x \Theta_n (\varphi_n(y, z) , z) }$, proving that $ \varphi_n$ is $C^2$ with respect to $y$. 
The remaining claims follow from Assumption \ref{ass:change of variable}. \ep

\medskip 
\begin{Lemma}\label{lemma:convergence of integrals}
Let $(\Yc, \Bc_{\Yc}, \nu_{\Yc})$ and $(\Zc, \Bc_{\Zc}, \nu_{\Zc})$ be finite-dimensional normed vector spaces endowed with their Borel $\sigma$-algebras and probability measures. 
Let $ 1 \leq \alpha \leq \infty$. 
Let $w \in \mathbb{L}^1 (\mathrm{d}x \otimes \nu_{\Zc}) \cap  \mathbb{L}^\infty (\mathrm{d}x \otimes \nu_{\Zc} ) $.
Let $ (k_n)$ be a sequence such that $ k_n :  \R \times \Yc \times \Zc \rightarrow \R $ satisfies
$$ \vert k_n(x, y, z) \vert \leq b_n(x) a_n(y, z),$$ 
where $ a_n \in \mathbb{L}^{\alpha}(\nu_{\Yc} \otimes \nu_{\Zc}) $ and $b_n \in  \mathbb{L}^{1}(\mathrm{d}x) \cap \mathbb{L}^{\infty}(\mathrm{d}x)$ are both positive. 
Finally, assume that $ (\Vert b_n\Vert_{\mathbb{L}^1 (\mathrm{d}x) } )_n$ and $(\Vert b_n\Vert_{\mathbb{L}^\infty (\mathrm{d}x) } )_n$ are bounded sequences and that, $(a_n)$ is $\alpha$-uniformly integrable in the sense that $ \lim_{ M \rightarrow \infty} \sup_n \Vert a_n \mathds{1}_{a_n \geq M} \Vert_{\mathbb{L}^{\alpha} (\nu_{\Yc} \otimes \nu_{\Zc} )} = 0$. 

\noindent Let $u_n: \R \times \Yc \times \Zc  \rightarrow \R$ satisfy Assumptions \ref{ass:change of variable} for $\alpha $ (and for the product measure for $w$) and $v_n : \R \times \Zc  \rightarrow \R$ satisfy Assumptions \ref{ass:change of variable} for $\alpha = 1$. 
Let $\varphi_n : \R \times \Yc \times \Zc \rightarrow \R$ (resp $\psi_n : \R \times \Yc \times \Zc  \rightarrow \R $) be defined by Lemma \ref{lemma:estim diffeo and change of variable} for the sequence $(u_n)$ (resp $(v_n)$). 
Then, using Notations \eqref{eqdef:convergencerandomvar}, we have the following expansion 
$$ 
\Big\Vert \int_{\Zc} \int_\R  k_n(x, \cdot , z) \big( w_n \big(x + u_n (x, \cdot, z ) , z\big) - w(x, z) \big)  \nu_{\Zc}( \mathrm{d} z ) \mathrm{d} x  \Big\Vert_{\mathbb{L}^{\alpha}(\nu_{\Yc})} \xrightarrow[n \rightarrow + \infty]{} 0
,$$
where $
w_n (x, z ) := w ( \psi_n( x , z ), z )  \partial_x \psi_n (x, z)
$.
\end{Lemma}

\noindent \textbf{Proof.} Let $M > 0$, since $w \in \mathbb{L}^{\infty}$, 
\begin{equation}\label{ineq:estim1 integral}
\Delta_n := \Big\vert 
\int_{\Zc} \int_\R  k_n(x, y , z) \big(  w_n ( x + u_n (x, y, z ) , z) - w(x, z) \big)   \mathrm{d} x  \, \nu_{\Zc}( \mathrm{d} z ) \Big\vert 
\leq 
I_1 (y) + I_2 (y)
\end{equation}
where $ I_1 := 2 \Vert w \Vert_{\mathbb{L}^{\infty}} \Vert b_n \Vert_{\mathbb{L}^1}
\int_{\Zc}  a_n(y, z) \mathds{1}_{\{a_n(y, z) \geq M \}} \nu_{\Zc}( \mathrm{d} z )$ and 
$$ I_2 := M \Vert b_n \Vert_{\mathbb{L}^\infty}
\int_{\Zc} \int_\R  \big\vert w_n \big(x + u_n (x, y, z ) , z\big) - w(x, z) \big\vert  \mathrm{d} x  \, \nu_{\Zc}( \mathrm{d} z ).$$ 
We first control $I_1$.
We have 
\begin{equation}
\Vert I_1 \Vert_{\mathbb{L}^{\alpha} (\nu_{\Yc})} \leq \Vert a_n \mathds{1}_{ {a_n \geq M}} \Vert_{\mathbb{L}^{\alpha} (\nu_{\Yc} \otimes \nu_{\Zc})},
\end{equation}
as a consequence of Jensen's inequality if $ 1 \leq \alpha < \infty $ and obvious if $ \alpha = +\infty$. 
Let $\varepsilon >0$. 
Now, we move on to bounding $I_2$. 
By assumption, $w \in \mathbb{L}^1 ( \nu_{\Zc} \otimes \mathrm{d}x )$, hence, by convolution, there exists a family $(w^\varepsilon)_\varepsilon$ such that 
\begin{equation}\label{ineq:estim w epsilon}
\int_{\Zc} \int_\R \vert w^\varepsilon (x, z ) - w (x, z) \vert \mathrm{d}x \,\nu_{\Zc} (\mathrm{d} z) \leq \varepsilon 
\,\, 
\text{and}
\,\, 
\vert w^\varepsilon (x, z) -  w^\varepsilon (y, z) \vert 
\leq 
C_\varepsilon \vert x - y \vert 
.\end{equation}
Set $ w_n^{\varepsilon} (x, z) := w^\varepsilon \big( \psi_n (x, z), z \big) \partial_x \psi_n (x, z) $. 
By the triangle inequality, 
\begin{equation}\label{ineq:estim2 integral}
I_2
\leq 
M \Vert b_n \Vert_{\mathbb{L}^\infty} \sum_{i = 1}^4  R_i (y)
\end{equation}
where 
\begin{equation*}
\begin{split}
R_1 (y) &:= 
\int_{\Zc} \int_\R \vert w^{\varepsilon}_n - w_n \vert (x + u_n (x, y, z ) , z)  \mathrm{d} x \nu_{\Zc}( \mathrm{d} z ) 
\\
R_2 (y) &:= 
\int_{\Zc} \int_\R \vert w^{\varepsilon}_n - w_n \vert (x , z) \mathrm{d} x \nu_{\Zc}( \mathrm{d} z ) \\
R_3 (y) &:=
\int_{\Zc} \int_\R  \vert w^{\varepsilon}_n (x + u_n (x, y, z ) , z)  - w^{\varepsilon}_n (x , z)  \vert \mathrm{d} x  \nu_{\Zc}( \mathrm{d} z ) 
\\
R_4 (y) &:= 
\int_{\Zc} \int_\R  \vert w^{\varepsilon}_n - w^{\varepsilon} \vert (x , z ) \mathrm{d} x \, \nu_{\Zc}( \mathrm{d} z ).
\end{split}
\end{equation*} 

\noindent By a change of variables and by Lemma \ref{lemma:estim diffeo and change of variable}, $\vert \partial_y \psi_n \vert \leq C$ and $\vert \partial_y \varphi_n \vert \leq C$, hence, using Estimate \eqref{ineq:estim w epsilon}, for $ i = 1, 2$,
\begin{equation}\label{ineq:estimate R2 R3}
R_i \leq C \int_{\Zc} \int_\R \vert w^{\varepsilon} - w \vert (x, z) \mathrm{d}x \, \nu_{\Zc} (\mathrm{d}z) \leq C \varepsilon 
.
\end{equation}
Furthermore, since $w^\varepsilon$ is Lipschitz in $x$, and $ \partial_x \psi_n $ is bounded, we have the existence of a constant $\hat{C}_\varepsilon >0$ such that 
\begin{equation}\label{ineq:estim R3}
R_3 \leq \hat{C}_\varepsilon \int_{\Zc} \int_\R  \vert u_n (x, y, z ) \vert \mathrm{d} x \, \nu_{\Zc}( \mathrm{d} z )  ,
\end{equation}
and, we also have (since one can choose $w^\varepsilon$ to be essentially bounded)
\begin{equation}\label{ineq:estim R4}
R_4 
\leq 
\Vert w^\varepsilon \Vert_{\mathbb{L}^\infty} \int_{\Zc} \int_\R  \vert \partial_x \psi_n (x, z ) - 1 \vert \mathrm{d} x \, \nu_{\Zc}( \mathrm{d} z )  
+
\hat{C}_\varepsilon \int_{\Zc} \int_\R  \vert \psi_n (x, z ) - x \vert \mathrm{d} x \, \nu_{\Zc}( \mathrm{d} z ) 
.\end{equation}
Now, putting all Estimates \eqref{ineq:estim R4}, \eqref{ineq:estim R3}, \eqref{ineq:estimate R2 R3}, \eqref{ineq:estim2 integral} and \eqref{ineq:estim1 integral} together, and using that $\Vert b_n \Vert_{\mathbb{L}^1 (\mathrm{d}x ) } $ and $\Vert b_n \Vert_{\mathbb{L}^\infty (\mathrm{d}x ) } $ are bounded, we get
\begin{align}
\Vert \Delta_n \Vert_{\mathbb{L}^{\alpha} (\nu_{\Yc} )}
\leq \,\, & M (\sup_n \Vert  b_n \Vert_{\mathbb{L}^\infty} ) \Big\{ \Vert w \Vert_{\mathbb{L}^\infty} \Vert \sup_{x \in \R} \{ \vert \partial_x \psi_n (x, \cdot) - 1\vert \} \Vert_{\mathbb{L}^{1}(\nu_{\Zc } )} \nonumber
\\
& \hspace{2.0cm}
+ \hat{C}_\varepsilon \Vert \sup_{x \in \R} \{ \vert \psi_n (x, \cdot) - x\vert \} \Vert_{\mathbb{L}^{1}(\nu_{\Zc}) } \nonumber
\\
& \hspace{2.0cm} + \hat{C}_\varepsilon \Vert \sup_{x \in \R} \{ \vert u_n (x, \cdot) - x\vert \} \Vert_{\mathbb{L}^{\alpha}( \nu_{\Yc} \otimes \nu_{\Zc} ) } + 2C \varepsilon \Big\} \nonumber
\\
&
+ \sup_n \Vert a_n \mathds{1}_{ \{a_n \geq M\}} \Vert_{\mathbb{L}^{\alpha} (\nu_{\Yc} \otimes \nu_{\Zc})} \label{eqdef:estimDletanLalpha}
.\end{align}
By Lemma \ref{lemma:estim diffeo and change of variable}, the first two terms in braces vanish as $n \rightarrow \infty$.
Taking the $\limsup$ in $n$ in \eqref{eqdef:estimDletanLalpha}, we get
$$\limsup_n \Vert \Delta_n \Vert_{\mathbb{L}^{\alpha} (\nu_{\Yc} )} 
\leq 
 M(\sup_n \Vert  b_n \Vert_{\mathbb{L}^\infty} ) 2C \varepsilon
+
M (\sup_n  \Vert a_n \mathds{1}_{ {a_n \geq M}} \Vert_{\mathbb{L}^{\alpha} (\nu_{\Yc} \otimes \nu_{\Zc})}) + 2 M \Vert a_n \Vert_{\mathbb{L}^\infty} C \varepsilon.$$
Now, letting $\varepsilon$ tend to $0$, and finally, letting $M$ tend to infinity gives the desired result. 
\ep

\medskip

\begin{Lemma}\label{lemma:expansion in Lp cdf}
 Let $\nu \in \Pi(\nu_1, \nu_2)$ for some $ \nu_1, \nu_2 \in \Pc_p (\R)$. Assume that $\nu$ admits the disintegration $\nu (\mathrm{d}x):= \nu_1 (\mathrm{d}x_1) k (x_1, x_2) \mathrm{d}x_2$ where $ k \in \mathbb{L}^\infty (\R^2)$. 
Let $ r_n \rightarrow 0 $, and let a family $(\Theta_n )$ be such that $ \Theta_n : \R^2 \rightarrow \R$ satisfies the following assumptions:
\begin{enumerate}[label=\textnormal{(\roman*)}]
\item For all $x_2 \in \R$, the mapping $ \Theta_n( \cdot, x_2 )$ is measurable. 
\item For all $ x_1 \in \R$, the mapping $ \Theta_n ( x_1, \cdot) $ is $C^2$ and, $ ( \vert \Theta_n \vert )_n $ is $p$-uniformly integrable. 
Furthermore,  $\text{max} \big( \vert \partial_{x_2} \Theta_n (x_1, x_2 ) \vert,  \vert \partial_{x_2 x_2}^2 \Theta_n (x_1, x_2 ) \vert \big) \leq C$.
\end{enumerate}
Define for $ n \in \N$  the following random variable, 
$$
Z_n
:=
F^{r_n} (X_2 + r_n \Theta_n ) ,
$$
where $F^{r_n} (x):= \nu \big(X_2 + r_n \Theta_n \leq x \big) $. 
Then we have the following expansion 
\begin{align*}
Z_n
&=
F_2(X_2) 
+
 r_n
 \int_\R
 (\Theta_n - \Theta_n (x_1, X_2) ) k(x_1, X_2) \nu_1(\mathrm{d} x_1)
+
 o_{\mathbb{L}^p (\nu)} (r_n) 
.\end{align*}
\end{Lemma}

 \noindent\textbf{Proof.} Let $r_n \rightarrow 0$ and $(\Theta_n)_n$ be as in the statement. 
 Set $\hat{\Theta}_n := r_n \Theta_n $, by Lemma \ref{lemma:estim diffeo and change of variable}, there exists $N \in \N$ and a family of functions $(\varphi_n)_n$ such that for $n \geq N$, $\varphi_n (x_1, \cdot)$ is the inverse of $x_2 \mapsto x_2 + \hat{\Theta}_n(x_1, x_2)$ and $\varphi_n$ is $C^2$ with respect to $x_2$ (here, as for $\hat{\Theta}_n$, the first argument of $\varphi_n$ is the parameter $x_1$, \textit{i.e.}, $\varphi_n (x_1, x_2)$ is the map of Lemma \ref{lemma:estim diffeo and change of variable} evaluated at $(y, z) = (x_2, x_1)$). 
 Set for $ n \geq N$, $k_n (x):= k \big(x_1, \varphi_n (x_1, x_2 )\big)  \partial_{x_2} \varphi_n ( x_1, x_2)  $ which is the disintegration with respect to the first marginal of the measure $ \nu \circ \big(X_1, X_2 + \hat{\Theta}_n \big)^{-1}$. 
 Finally, set
\begin{equation}\label{eqdef:decomposition difference}
\begin{split}
\Delta_n
&:=
Z_{n}
-
F_2(X_2) 
-r_n
\int_\R \big( \Theta_n (X) - \Theta_n (x_1, X_2) \big) k(x_1, X_2) \nu_1 (\mathrm{d}x_1)
,\end{split}
\end{equation}
By a change of variable $Z_{n} = \int_\R \int_{- \infty}^{X_2 + \hat{\Theta}_n (X)} k_n (x_1, x_2) \mathrm{d}x_2 \nu_1 (\mathrm{d}x_1)$.
Hence, we have the following equality 
\begin{equation} \label{eqdef:decomposition Z}
\begin{split}
Z_{n}
&=
 \int_\R \int_{- \infty}^{X_2} k_n (x) \mathrm{d}x_2 \nu_1 (\mathrm{d}x_1)
 +
 \int_\R \int_{X_2}^{X_2 + \hat{\Theta}_n (X)} k_n (x ) \mathrm{d}x_2 \nu_1(\mathrm{d}x_1)
 \\
 &=
  \int_\R \int_{- \infty}^{ \varphi_n (x_1, X_2 )} k(x) \mathrm{d}x_2 \nu_1(\mathrm{d}x_1)
 +
 \int_\R \int_{X_2}^{X_2 + \hat{\Theta}_n (X)} k_n (x ) \mathrm{d}x_2 \nu_1 (\mathrm{d}x_1)
. \end{split}
\end{equation}
Substituting $ Z_{n} $ by the expression \eqref{eqdef:decomposition Z} in the expression of $\Delta_{n} $ \eqref{eqdef:decomposition difference}, we obtain the following decomposition
\begin{equation*} 
\begin{split}
\Delta_{n}
&=
 \int_\R \int_{ X_2}^{ \varphi_n (x_1, X_2 )} k (x) \mathrm{d}x_2 \nu_1 (\mathrm{d}x_1) 
 +
 \int_\R \int_{X_2}^{X_2 + \hat{\Theta}_n(X)} k_n (x ) \mathrm{d}x_2 \nu_1 (\mathrm{d}x_1) 
 \\
 &-
  \int_\R \int_{X_2}^{X_2 + \hat{\Theta}_n(X)} k (x_1, X_2) \mathrm{d}x_2 \nu_1 (\mathrm{d}x_1) 
  -
  \int_\R \int_{ X_2}^{X_2 - \hat{\Theta}_n (x_1, X_2)} k (x_1, X_2) \mathrm{d}x_2 \nu_1 (\mathrm{d}x_1) 
  \\
  &=
   \int_\R \int_{ X_2 - \hat{\Theta}_n (x_1, X_2) }^{ \varphi_n(x_1, X_2 )} k (x) \mathrm{d}x_2 \nu_1 (\mathrm{d}x_1) 
    +
    \int_\R \int_{X_2}^{X_2 + \hat{\Theta}_n (X) } (k_n (x) - k (x_1, X_2) ) \mathrm{d}x_2 \nu_1 (\mathrm{d}x_1) 
  \\
 &+
  \int_\R \int_{ X_2}^{X_2 - \hat{\Theta}_n (x_1, X_2)} \big(k (x) - k(x_1, X_2) \big) \mathrm{d}x_2 \nu_1 (\mathrm{d}x_1) 
  = \sum_{i =1}^4 R_{i, n} , 
   \end{split}
\end{equation*}
where 
\begin{align*}
R_{1, n} &:=  \int_\R \int_{X_2 - \hat{\Theta}_n(x_1, X_2)}^{ \varphi_n ( x_1, X_2 ) }  k (x) \mathrm{d}x_2 \nu_1 (\mathrm{d}x_1) 
\,\, , \,\,
R_{2, n} := \hat{\Theta}_n (X) \int_\R (k_n - k) \big(x_1, X_2 \big) \nu_1 (\mathrm{d}x_1) 
\\
R_{3, n} &:= 
\int_\R \int_{X_2}^{X_2 + \hat{\Theta}_n(X) } \big(k_n (x) - k_n (x_1, X_2) \big) \mathrm{d}x_2 \nu_1 (\mathrm{d}x_1) 
\\
R_{4, n} &:=
 \int_\R \int_{X_2}^{X_2 - \hat{\Theta}_n (x_1, X_2)} \big(k (x) - k(x_1, X_2) \big) \mathrm{d}x_2 \nu_1 (\mathrm{d}x_1) .
\end{align*}

\smallskip

\noindent We will now prove that $ R_{i, n} = o_{\mathbb{L}^p (\nu ) } (r_n)$ for $ i = 1, \cdots, 4$. 
\textbf{For $i = 1$}, after a change of variable, 
\begin{align*}
R_{1, n} = \int_\R \int_0^1 
&k \big(x_1, (1 - \lambda) \big( X_2 - \hat{\Theta}_n (x_1, X_2) \big) + \lambda \varphi_n (x_1, X_2) \big) 
\\
&\times \big(\varphi_n (x_1, X_2) - X_2 + \hat{\Theta}_n (x_1, X_2) \big) \mathrm{d} \lambda \, \nu_1 (\mathrm{d}x_1)    
\end{align*}
and, by definition of $\varphi_n$, 
$$
\varphi_n (x_1, x_2) - x_2 = - \hat{\Theta}_n \big( x_1, \varphi_n (x_1, x_2) \big) 
=
- r_n \Theta_n \big( x_1, \varphi_n (x_1, x_2) \big) 
.$$
Hence, by the uniform Lipschitz property of $\Theta_n$ with respect to its second variable, and since, by Lemma \ref{lemma:estim diffeo and change of variable}, $ \Vert \sup_{x \in \R} \{ \vert \varphi_n (\cdot, x) - x \vert \} \Vert_{\mathbb{L}^p (\nu)} \rightarrow 0$, this yields the desired result.

\noindent \textbf{For $i = 2$}.  
After a change of variables, 
$$
R_{2, n} = r_n \Theta_n (X)\int_\R \big(k_n (x_1, X_2 \big) - k (x_1, X_2) \big) \nu_1 (\mathrm{d}x_1)
$$
which is a $ o_{\mathbb{L}^p (\nu)} (r_n) $ by an application of Lemma \ref{lemma:convergence of integrals} and since $\Theta_n$ is $p$-uniformly integrable.

\noindent \textbf{For $i = 3$}. 
After a change of variable, 
$$
R_{3, n} = \hat{\Theta}_n (X) 
\int_\R \int_{0}^{1} \big(k_n (x_1, (1 - \lambda)\hat{\Theta}_n(X) + X_2 \big) - k_n (x_1, X_2) \big) \mathrm{d}\lambda \nu_1 (\mathrm{d}x_1) 
:= R_{3, 1, n} + R_{3, 2, n}
,$$
where $ R_{3, 1, n} := \hat{\Theta}_n (X)\int_\R \int_{0}^{1} (k_n (x_1, (1 - \lambda)\hat{\Theta}_n(X) + X_2 ) - k (x_1, X_2) ) \mathrm{d}\lambda \nu_1 (\mathrm{d}x_1) $ and $ R_{3, 2, n} := \hat{\Theta}_n (X)\int_\R (k (x_1, X_2 ) - k_n (x_1, X_2) ) \nu_1 (\mathrm{d}x_1) $. By the case $i = 2$, we have $R_{3, 2, n} = o_{\mathbb{L}^p(\nu)}(r_n)$. 
Furthermore, since $k$ is bounded, $ \vert R_{3, 1, n}\vert \leq C \vert r_n \vert \vert \Theta_n \vert$. 
Hence, $ ( \vert \frac{R_{3, 1, n}}{r_n} \vert^p )_n$ is $p$-uniformly integrable.
Hence, we only need to check that $R_{3, 1, n} = o_{\mathbb{L}^1 (\nu ) } (r_n) $ and, 
$$
\frac{\E^{\nu} [ \vert R_{3, 1, n}\vert ]}{r_n}
   \leq   
 \int_{0}^1   \int_{\R}   \int_\R    \int_\R   \vert \Theta_n(z)
 k_n ( x_1, z_2 + \lambda \hat{\Theta}_n (z ) ) - k ( x_1, z_2) \vert \nu_1(\mathrm{d}x_1 ) \nu_1(\mathrm{d}z_1 ) \nu_2(\mathrm{d}z_2) \mathrm{d}\lambda.
$$
The right-hand side goes to $0$ by Lemma \ref{lemma:convergence of integrals}. 

\noindent \textbf{For $i = 4$}. 
By a change of variables, 
$$
R_{4, n} =
 \int_\R \int_{0}^1 \hat{\Theta}_n (x_1, X_2) \big(k (x_1, X_2) - k(x_1, X_2 + \lambda \hat{\Theta}_n (x_1, X_2)) \big)  \mathrm{d}\lambda \nu_1 (\mathrm{d}x_1).
$$
Hence, $\vert R_{4, n} \vert \leq r_n  \int_\R \int_{0}^1 \vert \Theta_n (x_1, X_2) \vert \vert k(x_1, X_2 + \lambda \hat{\Theta}_n (x_1, X_2)) - k (x_1, X_2) \vert  \mathrm{d}\lambda \, \nu_1 (\mathrm{d}x_1)$. Therefore, by Lemma \ref{lemma:convergence of integrals}, $\vert R_{4, n} \vert = o_{\mathbb{L}^p (\nu) }(r_n) $.
\ep

\appendix

\section{Discussion about Condition \hyperref[ass:mart coupling]{$\textbf{\rm C}_{\W_p^{\rm ad}}$} \ref{cond:info discrepancy}. } \label{subsec:discussion}

As in Subsection \ref{subsec:semi static hedging}, we identify $\mathbb{L}^1(\mu_i)$, $i = 1, 2$, with its image through the canonical injection $f \mapsto f(X_i)$ into $\mathbb{L}^1(\mu)$, which is an isometry since $\mu_i = \mu \circ X_i^{-1}$.
In the following, we identify the set of $\mu$-\textit{a.s.} constant functions with $\R$. Since $\mu, \mu_1$ and $\mu_2$ are probability measures, it is always true that 
$$
\R \subset \mathbb{L}^1(\mu_1) \cap \mathbb{L}^1(\mu_2)
.$$
We want to discuss for which additional condition on $\mu$ the converse inclusion (and the equality) holds.

\subsubsection{A toy example: product measure}

\begin{Proposition}\label{prop:independant measure}
Assume that $ \mu (\mathrm{d}x) = \mu_1(\mathrm{d}x_1) \otimes \mu_2(\mathrm{d}x_2) $ (where $\otimes$ denotes the product measure and $x = (x_1, x_2)$). 
In that case, we have 
$$\mathbb{L}^1(\mu_1) \cap \mathbb{L}^1(\mu_2) = \R.$$
\end{Proposition}

\proof
Let $ \mu (\mathrm{d}x) = \mu_1(\mathrm{d}x_1) \otimes \mu_2(\mathrm{d}x_2) $. 
Let $X = (X_1, X_2) \sim \mu $, then we have $X_1$ and $X_2$ which are independent. Let $H \in \mathbb{L}^1(\mu_1) \cap \mathbb{L}^1(\mu_2)$, which means that there exists $f \in \mathbb{L}^1(\mu_1) $ and $ g \in \mathbb{L}^1(\mu_2)$ such that 
$$
H(X) = f(X_1) = g(X_2) \,\, \mu\text{-a.s.}
.$$
Taking the conditional expectation (well defined since all functions are $\mu$-integrable) with respect to $X_1$ yields 
$$
f(X_1) = \E^\mu [g(X_2) \vert X_1 ] = \E^\mu [ g(X_2) ]
$$
where the last equality is a consequence of the independence. 
In that case we do have $f(X_1) = \E^\mu [ g(X_2) ]$, $\mu_1$-\textit{a.s.}, hence $f$ is $\mu_1$-\textit{a.s.} constant and we can conclude.
\ep

\subsubsection{If $\mu$ has a strictly positive density with respect to the Lebesgue measure and full support}

From now on we denote by $\Lebone $ the Lebesgue measure over $\R$ and $\Lebtwo $ the Lebesgue measure over $\R^2$.

\begin{Proposition}\label{full_suport_density}
Assume that $ \mu (\mathrm{d}x) = q(x) \Lebtwo(\mathrm{d}x)$, where $q(x) > 0 $ \textit{a.e.} on $\R^2$ (to be understood in the sense of the Lebesgue measure). In that case we have 
$$\R = \mathbb{L}^1(\mu_1) \cap \mathbb{L}^1(\mu_2) .$$
\end{Proposition}

\proof 
Since $q > 0$ \textit{a.e.} on $\R^2$, the marginals $\mu_1$ and $\mu_2$ admit densities with respect to $\Lebone$ which are positive \textit{a.e.} on $\R$. Hence $\mu$ is equivalent to $\Lebtwo$, and thus to the product measure $\mu_1 \otimes \mu_2$. In particular, $f(X_1) = g(X_2)$ holds $\mu$-\textit{a.s.} if and only if it holds $\mu_1 \otimes \mu_2$-\textit{a.s.}, and the result is immediate by the same reasoning as in the proof of Proposition \ref{prop:independant measure}.
\ep

\subsubsection{Positive density and support with connected interior}

\begin{Proposition}\label{connected}
Let $\Omega \subset \R^2$ be a Borel set with connected interior such that $\Lebtwo(\partial \Omega) =0$ (where $\partial$ denotes the boundary), and assume that $ \mu (\mathrm{d}x) = q(x) \mathds{1}_{\Omega} (x) \Lebtwo(\mathrm{d}x)$, where $q(x) > 0 $ \textit{a.e} on $\Omega$.
In that case we have 
$$\R = \mathbb{L}^1(\mu_1) \cap \mathbb{L}^1(\mu_2) .$$
\end{Proposition}

\proof 
Following the step of the previous proof, we get that Lebesgue almost everywhere on $\Omega$:
$$f(x_1) = g(x_2) .$$
Since we work up to Lebesgue-negligible sets, we may (and do) replace
$\Omega$ by its interior, which we still denote by $\Omega$; then
$\Omega$ is open and connected.
Now letting $I_1, I_2 \subset \R $, open and such that $ I_1 \times I_2 \subset \Omega$, we have that for almost all $x$ in $I_1 \times I_2$,
$f(x_1) = g(x_2)$.

In other words, letting
$$N = \{ (x_1, x_2) \in I_1 \times I_2 \,\,  \text{such that} \,\, f(x_1) \neq g(x_2) \,\, \}  \subset \R^2, $$
we have $ \Lebtwo(N) = 0.$
Again, following the steps of the last proofs, define
$$
N_{x_1} := \{ x_2 \in I_2 \,\, \text{such that} \,\, f(x_1) \neq g(x_2) \,\, \} $$
By Fubini's theorem, $\Lebtwo(N) = \int_{I_1} \Lebone(N_{x_1}) \Lebone(\mathrm{d}x_1) = 0$.
Hence, letting
$$
A_1 = \{ x_1 \in I_1 \,\, \Lebone(N_{x_1}) = 0 \,\, \},
$$
we have $\Lebone(I_1 \setminus A_1) = 0$.
Let $x_1 \in A_1$, we have that for almost all $x_2 \in I_2$, $f(x_1) = g(x_2)$, hence $g$ is almost everywhere constant on $I_2$.
By similar considerations, $f$ is almost everywhere constant on $I_1$.
Hence, there exists a constant $c_{I_1 \times I_2}$ such that
$$
f(x_1) = g(x_2) = c_{I_1 \times I_2} \,\, \text{$\Lebtwo$  a.e on $I_1 \times I_2$}
.$$

\medskip

Now let $I_1, I_2 \subset \R $ and $\hat{I}_1, \hat{I}_2 \subset \R$ such that $I_1\times I_2 \subset \Omega$ and $\hat{I}_1 \times \hat{I}_2 \subset \Omega $ satisfy $ \Lebtwo \big( (I_1 \times I_2) \cap (\hat{I}_1 \times \hat{I}_2 ) \big) >0$.
Then picking $x$ in $\big( (I_1 \times I_2) \cap (\hat{I}_1 \times \hat{I}_2 ) \big)$ (which exists since this set has strictly positive measure), the constant must coincide and
\begin{equation}\label{localconstance}
c_{I_1 \times I_2} = f(x_1) = g(x_2) = c_{\hat{I}_1 \times \hat{I}_2}
.
\end{equation}
Hence, the function that maps $ \Gamma : \Omega \rightarrow \R $, mapping $x \in \Omega$ to $c_{I_1 \times I_2}$ where $I_1\times I_2$ is a neighborhood of $x$ in $\Omega$ is well defined because $\Omega$ is open and the constant does not depend on the neighborhood $I_1\times I_2$.
Indeed, consider another neighborhood $\hat{I}_1 \times \hat{I}_2$, the intersection $(I_1\times I_2) \cap ( \hat{I}_1 \times \hat{I}_2)$ remains an open neighborhood of $x$, which will have non zero $\Lebtwo$ measure. Hence, $\Gamma(x) = c_{I_1 \times I_2} = c_{\hat{I}_1 \times \hat{I}_2}$ by equation \eqref{localconstance}.
And furthermore, $\Gamma$ is locally constant, since for $x \in \Omega$, one can find a neighborhood $I_1 \times I_2$ of $x$ where for all $y \in \Omega$, $\Gamma(x) = \Gamma(y) = c_{I_1 \times I_2}$.
Hence $\Gamma$ is a locally constant function on a connected set $\Omega$, so $\Gamma$ is constant and $c$ does not depend on $I_1 \times I_2$.
\ep

\subsubsection{A general condition}

In this subsection, we let $\Omega := \text{supp} (\mu)$ and $\Omega_1 := \text{supp} (\mu_1)$, and, for $x_1 \in \Omega_1$, let $I_{x_1} := \text{supp} (\pi_{x_1})$ where  $(\pi_{x_1})_{x_1 \in \Omega_1}$ denotes the disintegration of $\mu$ with respect to $\mu_1$. 

\begin{Proposition}\label{gen_prop}
There exists $\hat{\Omega} \subset \Omega $ such that $\mu ( \Omega \setminus \hat{\Omega}) = 0$ and such that for all $x,x'\in \hat{\Omega}$, there exists
$N\in\mathbb{N}$ and points $x=x^1,\dots,x^N = x'\in \hat{\Omega}$ such that
$$
\mu_2\bigl( I_{x_1^i} \cap I_{x_1^{i+1}} \bigr) > 0,
\qquad\text{for } 1\le i \le N-1 .
$$
In that case we have 
$$\R = \mathbb{L}^1(\mu_1) \cap \mathbb{L}^1(\mu_2) .$$
\end{Proposition}

\begin{Remark}
We understand that the property $\R = \mathbb{L}^1(\mu_1) \cap \mathbb{L}^1(\mu_2)$ is in fact really dependent on the support of $\mu$. Let $
\mu = \mathcal{L}(X_1, X_1 + Z),
$ where $
\mathcal{L}(X_1)(\mathrm{d}x_1)
  = \frac{1}{2}\bigl(\delta_{-1}(\mathrm{d}x_1)
                     + \delta_{1}(\mathrm{d}x_1)\bigr),
\qquad
\mathcal{L}(Z) = \mathcal{U}([-1,1])
$. 
The support of $\mu$ is given by 
\begin{center}
\begin{tikzpicture}[scale=1.0]
  \draw[->] (-2.5,0) -- (2.5,0) node[below left] {$X_1$};
  \draw[->] (0,-2.5) -- (0,2.5) node[above right] {$X_1 + Z$};

  \foreach \x in {-1,1}{
    \draw (\x,0.07) -- (\x,-0.07);
    \node[below] at (\x,0) {\small $\x$};
  }

  \foreach \y in {-2,-1,1,2}{
    \draw (0.07,\y) -- (-0.07,\y);
    \node[left] at (0,\y) {\small $\y$};
  }

  \draw[very thick,blue] (1,0) -- (1,2);
  \node[blue,right] at (1,1) {$X_1 = 1$};

  \draw[very thick,red] (-1,-2) -- (-1,0);
  \node[red,left] at (-1,-1) {$X_1 = -1$};

  \foreach \p in {(1,0),(1,2),(-1,-2),(-1,0)}{
    \fill \p circle (1.3pt);
  }
\end{tikzpicture}
\end{center}
Or similarly, consider the measure $\mu = \Lc ( X, X + ZX) $ where $X,Z \sim \mathcal{U}([-1, 1])$ are independent.
\begin{center}
    \begin{tikzpicture}[scale=1.4]
  \draw[->] (-1.5,0) -- (1.5,0) node[below left] {$X_1$};
  \draw[->] (0,-2.5) -- (0,2.5) node[above right] {$X_1 + X_1 Z$};

  \foreach \x in {-1,1}{
    \draw (\x,0.07) -- (\x,-0.07);
    \node[below] at (\x,0) {\small $\x$};
  }

  \foreach \y in {-2,-1,1,2}{
    \draw (0.07,\y) -- (-0.07,\y);
    \node[left] at (0,\y) {\small $\y$};
  }

  \fill[blue!15] (0,0) -- (1,0) -- (1,2) -- cycle;
  \fill[blue!15] (0,0) -- (-1,0) -- (-1,-2) -- cycle;

  \foreach \p in {(-1,0),(1,0),(-1,-2),(1,2),(0,0)}{
    \fill \p circle (1.2pt);
  }
\end{tikzpicture}
\end{center}
In that case again we have $\text{sgn}(X_1) = \text{sgn}(X_2)$. We see that this measure does not fall under the scope of Proposition \ref{gen_prop} because the support of the conditional expectation collapses at $x_1 = 0$, which separates the support into two connected components. 
We see that all cases proven in Propositions \ref{prop:independant measure}, \ref{full_suport_density} and \ref{connected} fall under the scope of Proposition \ref{gen_prop}. 
Indeed, under the assumptions of Propositions \ref{prop:independant measure}, \ref{full_suport_density}, we have that for all $x_1 \in \Omega_1$, $\text{supp}( \Lc(X_2 \vert X_1) ) = \text{supp}(\Lc(X_2))$. 
In the case of Proposition \ref{connected}, we have $\mu( \Omega \setminus \text{Int}(\Omega) )=0$.
Since $\text{Int}(\Omega)$ is connected, for all $x, x' \in \text{Int}(\Omega)$,
there exists a continuous path $\gamma : [0,1] \mapsto \text{Int}(\Omega)$ with $\gamma(0) = x$ and $\gamma(1) = x'$.
By continuity of $\gamma$ and compactness of $[0,1]$, one can cover the path by a finite number of balls. In other words, there exists $N \in \N$ and  $x = x^1 , \cdots, x^N = x'$, and $r_1, \cdots, r_N > 0$ such that for $1 \leq i \leq N$, 
$
B_i := B(x^i, r_i) \subset \text{Int} (\Omega) 
$ and for $1\leq i \leq N-1$, $ B_i \cap B_{i+1} \neq \emptyset$. In that case the points $x^1, \cdots, x^N$ satisfy the conditions of proposition \ref{gen_prop}.
\end{Remark}

\proof
Let $h \in \mathbb{L}^1(\mu_1) \cap \mathbb{L}^1(\mu_2)$.  
Then there exist $f\in \mathbb{L}^1(\mu_1)$ and $g\in \mathbb{L}^1(\mu_2)$ such that
$$
h(x_1,x_2)=f(x_1)=g(x_2)
\qquad \mu\text{-a.e. on }\Omega.
$$
Let
$$
N:=\{(x_1,x_2)\in\Omega : f(x_1)\neq g(x_2)\}.
$$
Since $\mu(N)=0$, Fubini's theorem yields a set 
$A\subset\Omega_1$ with $\mu_1(A)=0$ such that for all 
$x_1\in\Omega_1\setminus A$,
$$
\pi_{x_1}(N_{x_1})=0,
\,\,
N_{x_1} := \{x_2 : (x_1,x_2)\in N\}.
$$
Hence, for every $x_1\notin A$, 
$$
f(x_1)=g(x_2) \,\, \pi_{x_1}-\textit{a.s.} \text{ on } I_{x_1}.
$$
In other words there exists $c_{x_1}$ such that $g$ is constant equal to $c_{x_1}$ $\pi_{x_1}$ almost surely.
It is also clear that if $x_1$ and $x_1'$ are such that $\mu_{2}(I_{x_1} \cap I_{x_1'}) >0$, in that case $c_{x_1} = c_{x_1'}$. 
Using the chain condition stated in Proposition \ref{gen_prop} yields the desired result. \ep

\normalem
\printbibliography

\end{document}